\documentclass{amsart}
\usepackage{amsmath,amsthm}
\usepackage{amsfonts,amssymb}

\newtheorem{thm}{Theorem}[section]
\newtheorem{cor}[thm]{Corollary}

\newtheorem{prop}[thm]{Proposition}

\newtheorem{defn}[thm]{Definition}

\theoremstyle{remark}

 \def\CD{{\mathcal D}}
 
 \def\CH{{\mathcal H}}

 \def\CP{{\mathcal P}}     
 \def\CM{{\mathcal M}}     
      
 \def\CV{{\mathcal V}}

 \def\NN{{\mathbb N}}

 \def\RR{{\mathbb R}}
 \def\ZZ{{\mathbb Z}}
        \def\co{\operatorname{co}}
        \def\proj{\operatorname{proj}}
        \def\supp{\operatorname{supp}}

\begin{document}
 
\title[almost everywhere convergence]
  {Almost Everywhere Convergence of orthogonal expansions of several variables}
\author{Yuan Xu}
\address{Department of Mathematics\\ University of Oregon\\
    Eugene, Oregon 97403-1222.}\email{yuan@math.uoregon.edu}

\date{May 28, 2003, revised Oct. 11, 2003.}
\thanks{Work supported in part by the National Science 
Foundation under Grant DMS-0201669}
                                                    
\begin{abstract}
For weighted $L^1$ space on the unit sphere of $\RR^{d+1}$, in which the 
weight functions are invariant under finite reflection groups, a maximal 
function is introduced and used to prove the almost everywhere convergence 
of orthogonal expansions in $h$-harmonics. The result applies to various
methods of summability, including the de la Vall\'ee Poussin means and the 
Ces\`aro means. Similar results are also established for weighted orthogonal 
expansions on the unit ball and on the simplex of $\RR^d$. 
\end{abstract}

\maketitle                      
 
\section{Introduction} 
\setcounter{equation}{0}

Let $S^d = \{x: \|x\| =1\}$ denote the unit sphere in $\RR^{d+1}$, where 
$\|x\|$ denote the usual Euclidean norm. For a nonzero vector $v \in 
\RR^{d+1}$, let $\sigma_v$ denote the reflection with respect to the 
hyperplane perpendicular to $v$, $x \sigma_v : = x - 2 (\langle x,v \rangle 
/\|v\|^2) v$, $x \in \RR^{d+1}$, where $\langle x,y\rangle$ denote the usual 
Euclidean inner product. We consider the weighted approximation on $S^d$ with 
respect to the measure $h_\kappa^2 d \omega$, where $d \omega$ is the surface 
(Lebesgue) measure on $S^d$ and the weight function $h_\kappa$ is defined by 
\begin{equation}\label{eq:1.1}
h_\kappa(x) = \prod_{v \in R_+} |\langle x, v\rangle|^{\kappa_v}, \qquad 
   x \in \RR^{d+1}, 
\end{equation}
in which $R_+$ is a fixed positive root system of $\RR^{d+1}$, normalized so 
that $\langle v, v \rangle =2$ for all $v \in R_+$, and $\kappa$ is a 
nonnegative multiplicity function $v \mapsto \kappa_v$ defined on $R_+$ with 
the property that $\kappa_u = \kappa_v$ whenever $\sigma_u$ is conjugate to 
$\sigma_v$ in the reflection group $G$ generated by the reflections 
$\{\sigma_v:v \in R_+\}$. The function $h_\kappa$ is invariant under the 
reflection group $G$. The simplest example is given by the case 
$G=\ZZ_2^{d+1}$ for which $h_\kappa$ is just the product weight function
\begin{equation}\label{eq:1.2}
  h_\kappa (x)  = \prod_{i=1}^{d+1} |x_i|^{\kappa_i}, \qquad \kappa_i \ge 0. 
\end{equation}    
We denote by $a_\kappa$ the normalization constant of $h_\kappa$, 
$a_\kappa^{-1}  = \int_{S^d} h_\kappa^2(y) d\omega$, and denote by 
$L^p(h_\kappa^2)$, $1\le p\le \infty$, the space of functions defined on 
$S^d$ with the finite norm
$$
\|f\|_{\kappa,p} := \Big(a_\kappa \int_{S^d} |f(y)|^p h_\kappa^2(y) d\omega(y) 
\Big)^{1/p}, \qquad  1 \le p < \infty,
$$
and for $p = \infty$ we assume that $L^\infty$ is replaced by $C(S^{d})$, the 
space of continuous functions on $S^d$ with the usual uniform norm 
$\|f\|_\infty$. We consider the weighted approximation in $L^p(h_\kappa^2)$.
The case $\kappa \equiv 0$ corresponds to the usual $L^p$ (unweighted) 
approximation on $S^d$.

The homogeneous polynomials that are orthogonal with respect to $h_\kappa^2
d\omega$ are studied by Dunkl (\cite{D1,D2}; see \cite{DX} and the references 
therein). They are called $h$-harmonics, since they satisfy many properties 
that are similar to those of ordinary harmonics. In particular, summability 
of $h$-harmonic expansions has been studied in \cite{X97b,LX,X01a} and
weighted approximation theory by polynomials in $L^p(h_\kappa^2)$ has been
developed in \cite{X02,X02b}. For the usual harmonic analysis and approximation
on the sphere, see \cite{BC,LW,LN,Rus} and the references therein. The study 
in the weighted case often becomes more difficult, since the 
orthogonal group acts transitively on the sphere $S^d$ but a reflection group 
does not.  This happens to the case of almost everywhere convergence, which we 
study in this paper. 

Let $d(x,y) = \arccos \langle x,y\rangle$ be the geodesic distance of 
$x,y \in S^d$. For the usual approximation on $S^d$, the maximal function 
\begin{equation}\label{eq:1.3}
\CM f(x):= \sup_{0\le \theta \le \pi}\frac{\int_{d(x,y)\le \theta}f(y)d\omega}
    {\int_{d(x,y)\le \theta}d\omega}  
\end{equation}
is a weak type $(1,1)$ operator and it plays an important role in the study
of almost everywhere convergence. As we will see, however, the straightforward 
extension 
$$
  \CM_h f(x) = \sup_{0\le \theta \le \pi} 
   \frac{\int_{d(x,y)\le \theta} f(y)h_\kappa^2(y)d\omega} 
          {\int_{d(x,y) \le \theta} h_\kappa^2(y)d\omega} 
$$
for the weighted $L^1$ functions is not the natural one for the weighted 
approximation. One of our main results is to introduce an alternative maximal 
function, making use of the weighted spherical means $T_\theta^\kappa$ studied
in \cite{X02b}, and show that it is weak type $(1,1)$ in a proper sense. The 
weak type inequality is then used to prove the almost everywhere convergence of
the $h$-harmonic expansions. In particular, as corollaries, our results show 
that the Poisson integral, the Ces\`aro means and the de la Vall\`ee
Poussin means of the $h$-harmonic expansions all converge almost everywhere 
on the sphere for $f \in L^1(h_\kappa^2)$. 

Our main result also includes extension of these results to weighted 
approximation on the unit ball $B^d = \{x: \|x\|\le 1\}$, $x \in \RR^d$, in 
which the weight function is of the form
\begin{equation}\label{eq:1.4}
 W_{\kappa,\mu}^B(x) = h_\kappa^2(x) (1-\|x\|^2)^{\mu-1/2}, \qquad x \in B^d,
\end{equation}
where $h_\kappa$ is a reflection invariant weight function on $\RR^d$ and
$\mu \ge 0$, and to the weighted approximation on the simplex 
$$
T^d = \{x \in \RR^d:x_1 \ge 0, \ldots, x_d \ge 0, 1-|x| \ge 0\}, \qquad
 |x| = x_1+\ldots + x_d,
$$
in which the weight functions take the form
\begin{equation}\label{eq:1.5}
W_{\kappa,\mu}^T(x) = h_\kappa^2(\sqrt{x_1}, \ldots,\sqrt{x_d})
   (1-|x|)^{\mu-1/2} /\sqrt{x_1 \cdots x_d},  
\end{equation}
where $\mu \ge 0$ and $h_\kappa$ is a reflection invariant weight function 
defined on $\RR^d$, even in each of its variables. These include the
classical weight functions on these domains, which are
\begin{equation}\label{eq:1.6}
W_\mu^B(x) = (1-\|x\|^2)^{\mu-1/2}, \qquad x \in B^d,
\end{equation}
on the unit ball (taking $h_\kappa (x) =1$) and 
\begin{equation}\label{eq:1.7}
W_\kappa^T(x) = x_1^{\kappa_1-1/2} \cdots x_d^{\kappa_d-1/2} 
     (1-|x|)^{\kappa_{d+1}-1/2},\qquad x \in T^d,
\end{equation} 
on the simplex (taking $h_\kappa(x) = \prod_{i=1}^d|x_i|^{2\kappa_i}$ and
$\kappa_{d+1} = \mu$). While the results on $B^d$ can be established easily
using the corresponding result on $S^d$, those on $T^d$ need more careful 
study. Our results on almost everywhere convergence appear to be new even for 
the classical weight functions. For $d =1$, the two cases correspond to the 
Gegenbauer expansions and to the Jacobi expansions, respectively. 

The paper is organized as follows. The maximal function and almost everywhere
convergence on the unit sphere are studied in Section 2. The results on the
unit ball and on the simplex are studied in Section 3 and Section 4, 
respectively.  

\section{Maximal function and almost everywhere convergence on $S^d$}
\setcounter{equation}{0}

\subsection{Background}
Let $h_\kappa$ be the reflection invariant weight function defined in 
\eqref{eq:1.1}. The essential ingredient of the theory of $h$-harmonics is 
a family of first-order differential-difference operators, $\CD_i$, called 
Dunkl's operators, which generates a commutative algebra; these operators are 
defined by (\cite{D1})
$$
  \CD_i f(x) = \partial_i f(x) + \sum_{v \in R_+} k_v 
    \frac{f(x) -  f(x \sigma_v)} {\langle x, v\rangle}
        \langle v,\varepsilon_i\rangle, \qquad 1 \le i \le d+1, 
$$
where $\varepsilon_1, \ldots, \varepsilon_d$ are the standard unit vectors of 
$\RR^{d+1}$. The $h$-Laplacian is defined by $\Delta_h=\CD_1^2 + \ldots + 
\CD_{d+1}^2$ and it plays the role similar to that of the ordinary Laplacian. 
Let $\CP_n^{d+1}$ denote the subspace of homogeneous polynomials of degree $n$ 
in $d+1$ variables. An $h$-harmonic polynomial $P$ of degree $n$ is a 
homogeneous polynomial $P \in \CP_n^{d+1}$ such that $\Delta_h P  =0$. 
Furthermore, let $\CH_n^{d+1}(h_\kappa^2)$ denote the space of $h$-harmonic 
polynomials of degree $n$ in $d+1$ variables and define 
$$
 \langle f, g\rangle_\kappa : = a_\kappa \int_{S^d} f(x) g(x) 
    h^2_\kappa(x) d\omega(x). 
$$ 
Then $\langle P,Q\rangle_{\kappa} = 0$ for $P \in \CH_n^{d+1}(h_\kappa^2)$ and 
$Q \in \Pi_{n-1}^{d+1}$. The spherical $h$-harmonics are the restriction of 
$h$-harmonics on the unit sphere. It is known that $\dim \CH_n^{d+1}
(h_\kappa^2) =\dim \CP_n^{d+1} -\dim \CP_{n-2}^{d+1}$ with $\dim \CP_n^d = 
\binom{n+d-1}{n}$. 

The standard Hilbert space theory shows that
$$
L^2(h_\kappa^2) = \sum_{n=0}^\infty\bigoplus \CH_n^{d}(h_\kappa^2).
$$
That is, with each $f\in L^2(h_\kappa^2)$ we can associate its $h$-harmonic
expansion
$$
  f(x) = \sum_{n=0}^\infty Y_n(h_\kappa^2;f,x), \qquad x \in S^d,
$$
in $L^2(h_\kappa^2)$ norm. For the surface measure ($\kappa =0$), such a
series is called the Laplace series (cf. \cite[Chapt. 12]{Er}). The orthogonal
projection $Y_n(h_\kappa^2): L^2(h_\kappa^2) \mapsto \CH_n^{d+1}(h_\kappa^2)$
takes the form
\begin{equation}\label{eq:2.1}
 Y_n(h_\kappa^2;f,x):=
   \int_{S^d} f(y) P_n(h_\kappa^2;x,y) h_\kappa^2(y)\, d\omega(y),
\end{equation}
where the kernel $P_n(h_\kappa^2;x,y)$ is the reproducing kernel of the space 
of $h$-harmonics $\CH_n^{d+1}(h_\kappa^2)$ in $L^2(h_\kappa^2)$. The kernel 
$P_n(h_\kappa^2;x,y)$ enjoys a compact formula in terms of the intertwining 
operator between the commutative algebra generated by the partial derivatives 
and the one generated by Dunkl's operators. This operator, $V_\kappa$, is 
linear and it is determined uniquely by
$$
 V_\kappa \CP_n^d \subset \CP_n^d, \qquad V_\kappa 1=1,
    \qquad \CD_i V_\kappa = V_\kappa \partial_i,  \qquad 1 \le i \le d+1.
$$
The compact formula of the reproducing kernel for $\CH_n^{d+1}(h_\kappa^2)$
is given by (\cite{X97b})
\begin{equation}\label{eq:2.2}
P_n(h_\kappa^2;x,y) = \frac{n+\lambda_\kappa}{\lambda_\kappa}
   V_\kappa[C_n^{\lambda_\kappa} (\langle \cdot, y \rangle )](x),
\end{equation}
where, and throughout this paper, we fix the value of $\lambda_\kappa$ as 
\begin{equation}\label{eq:2.3}
  \lambda := \lambda_\kappa = \gamma_\kappa+ \frac{d-1}{2} \qquad 
 \hbox{with} \qquad  \gamma_\kappa =  \sum_{v\in R_+} \kappa_v.
\end{equation}
The function $C_n^\lambda(t)$ is the standard Gegenbauer polynomial, 
orthogonal with respect to the weight function
$$
w_\lambda(t) = (1-t^2)^{\lambda-1/2}, \qquad -1\le t \le 1.
$$
If all $\kappa_v = 0$, $V_\kappa$ becomes the identity operator and the above
formula is the usual formula for the zonal polynomial.

An explicit formula of $V_\kappa$ is known only in the case of symmetric group
$S_3$ for three variables and in the case of the abelian group $\ZZ_2^{d+1}$. 
In the latter case, $V_\kappa$ is an integral operator,
\begin{equation} \label{eq:2.4}
  V_\kappa f(x) = c_\kappa
       \int_{[-1,1]^{d+1}} f(x_1 t_1, \ldots,x_{d+1} t_{d+1}) 
        \prod_{i=1}^{d+1} (1+t_i)(1-t_i^2)^{\kappa_i -1} d t,
\end{equation}
where $c_\kappa$ is the normalization constant determined by $V_\kappa 1 =1$.
If some $\kappa_i =0$, then the formula holds under the limit relation
$$
 \lim_{\lambda \to 0} c_\lambda \int_{-1}^1 f(t) (1-t^2)^{\lambda -1} dt
  = [f(1) + f(-1)] /2.
$$
One important property of the intertwining operator is that it is positive
(\cite{Ros}), that is, $V_\kappa p \ge 0$ if $p \ge 0$. In fact, for each 
$x \in \RR^{d+1}$ there is a unique probability measure $\mu_x$ such that 
\begin{equation} \label{eq:2.5}
     V_\kappa f(x) = \int_{\RR^{d+1}} f(y) d\mu_x(y)    
\end{equation}
for each polynomial $f$. The measure $\mu_x$ is compactly supported with 
$\supp \mu_x \subset \co\{w x: w\in G\}$, the convex hull of the orbit of
$x$ under $G$. The intertwining operator plays an essential role for the 
weighted approximation by polynomials in $L^p(h_\kappa^2)$, as shown in 
\cite{X02b}, and it is also essential for the definition of our maximal 
function. 

\subsection{Maximal Function} 
For the definition of our maximal function, we need the extension of the 
spherical means defined in \cite{X02}. These means are defined implicitly
as follows:

\begin{defn}\label{defn1}
For $0 \le \theta \le \pi$, the means $T_\theta^\kappa$ are defined by 
\begin{equation}\label{eq:defn1}
c_\lambda \int_0^\pi T_\theta^\kappa f(x) g(\cos \theta)
  (\sin \theta)^{2\lambda} d\theta
= a_\kappa \int_{S^d} f(y) V_\kappa [g(\langle x,\cdot \rangle)](y) 
   h_\kappa^2(y)d\omega(y),
\end{equation} 
where $g$ is any $L^1(w_\lambda)$ function and $\lambda = \lambda_\kappa$.
\end{defn} 

The weighted spherical means $T_\theta^\kappa$ are well-defined (\cite{X02b}). 
In the case of the Lebesgue measure (that is, $h_\kappa(x)=1$), $\kappa =0$ 
and $V_\kappa = id$, the means $T_\theta^\kappa$ agree with the usual spherical
means $T_\theta f$, 
\begin{equation*}  
 T_\theta f(x) = \frac{1}{\sigma_{d-1} (\sin \theta)^{d-1}} 
     \int_{\langle x, y\rangle = \cos\theta} f(y) d\omega(y), 
\end{equation*}  
where $\sigma_{d-1} = \int_{S^{d-1}} d\omega = 2 \pi^{d/2}/\Gamma(d/2)$ is 
the surface area of $S^{d-1}$. The properties of the means $T_\theta$ are 
well-known; see \cite{BBP,P}, for example. The means $T_\theta^\kappa$ satisfy
similar properties as shown in \cite{X02,X02b}. In particular, for 
$f\in L^p(h_\kappa^2)$, $1 \le p < \infty$, or $f\in C(S^d)$,
\begin{equation}\label{eq:2.7}
 \|T_\theta^\kappa f\|_{\kappa,p} \le \|f\|_{\kappa,p} \qquad \hbox{and}
  \qquad \lim_{\theta \to 0} \|T_\theta^\kappa f -f \|_{\kappa,p} =0.
\end{equation}
Let us also mention that if $f_0(x) =1$, then $T_\theta^\kappa f_0(x) =1$. In 
particular, setting $f(x) =1$ in \eqref{eq:defn1} gives the equation
\begin{equation}\label{eq:2.8}
 a_\kappa \int_{S^d} V_\kappa [g(\langle x,\cdot \rangle)](y) 
   h_\kappa^2(y)d\omega(y) = c_\lambda \int_0^\pi g(\cos \theta)
  (\sin \theta)^{2\lambda} d\theta
\end{equation} 
for every $g \in L^1(w_\lambda)$, which is a special case proved early in
\cite{X97b}. A more general result proved in \cite{X97b} is the following
formula
\begin{equation}\label{eq:2.8.5}
 a_\kappa \int_{S^d} V_\kappa f(y) h_\kappa^2(y) d\omega =
    A_\kappa \int_{B^{d+1}} f(x) (1-\|x\|^2)^{\gamma_\kappa -1} dx, 
\end{equation}
where $A_\kappa$ is the normalization constant for 
$(1-\|x\|^2)^{\gamma_\kappa-1}$. The spherical means are used to define a 
modulus of continuity of $f$, which is used in turn to characterize the best 
approximation for polynomials. We use them to define our maximal function.

\begin{defn}\label{defn2}
For $f\in L^1(h_\kappa^2)$, the maximal function $\CM_\kappa f$ is defined by 
\begin{equation}\label{eq:defn2}
  \CM_\kappa f(x) = \sup_{0\le \theta \le \pi} 
    \frac{\int_0^\theta T_\phi^\kappa |f|(x) (\sin \phi)^{2\lambda} d\phi}
      {\int_0^\theta (\sin \phi)^{2\lambda} d\phi}.
\end{equation} 
\end{defn} 

Let $c(x,\theta):= \{y \in B^{d+1}: \langle x,y\rangle \ge \cos \theta\}$.
Its restriction on $S^d$, which we denoted by $c_S(x,\theta)$, is the 
spherical cap centered at $x$. Let $\chi_E$ denote the 
characteristic function of the set $E$. For $0 \le \theta \le \pi$, let 
$g_\theta(t) = \chi_{[\cos\theta,1]}(t)$ be a function defined on $[-1,1]$. 
It follows from the definition of $T_\theta^\kappa$ that 
$$
c_\lambda \int_0^\theta T_\phi^\kappa |f|(x) (\sin \phi)^{2\lambda} d\phi
= a_\kappa \int_{S^d} |f(y)| V_\kappa [g_\theta(\langle x,\cdot \rangle)](y) 
   h_\kappa^2(y)d\omega(y).
$$
On the other hand, it is easy to see that $g_\theta(\langle x,\cdot \rangle)=
\chi_{c(x,\theta)}$. Hence, using the equation 
\eqref{eq:2.8}, we have the following alternative definition of 
$\CM_\kappa f(x)$. 

\begin{prop} \label{prop:2.3}
For $f\in L^1(h_\kappa^2)$, 
$$
  \CM_\kappa f(x) =  \sup_{0\le \theta \le \pi} \frac{
  \int_{S^d} |f(y)| V_\kappa \left[\chi_{c(x,\theta)}\right](y)
 h_\kappa^2(y)d\omega(y)}
      {\int_{S^d} V_\kappa \left[\chi_{c(x,\theta)}\right](y)
h_\kappa^2(y)d\omega(y)}
$$
\end{prop}

In the case of $\kappa =0$, $V_\kappa = id$; we see that $\CM_\kappa f$ 
reduces to the classical maximal function \eqref{eq:1.3} for $|f|$. In the 
case of $G = \ZZ_2^{d+1}$, the formula of $V_\kappa$ in \eqref{eq:2.4}
shows that 
$$
V_\kappa \left[\chi_{c(x,\theta)}\right](y) = 0, \qquad \hbox{if} \qquad 
   \langle \bar x, \bar y \rangle < \cos \theta, 
$$
where $\bar x = (|x_1|, \ldots,|x_d|)$. Evidently, $\CM_\kappa f(x)$ satisfies 
$$
  \|\CM_\kappa f \|_\infty \le \|f\|_\infty, \qquad\quad \forall f \in C(S^d).
$$
Our goal is to prove a weak type estimate for $f \in L^1(h_\kappa^2)$. In the 
following, the constant $c$ denotes a generic constant, which depends only on 
the values of $d$, $\kappa$ and other fixed parameters and whose value may be 
different from line to line. 

\begin{thm} \label{thm:maximal}
Let $\sigma > 0$ and $f\in L^1(h_\kappa^2)$. Define $c(\sigma) = 
\{x: \CM_\kappa f(x) \ge \sigma\}$. Then there is a function $w_\sigma(x)$ 
which is positive on $c(\sigma)$ such that 
$$
\int_{S^d} \chi_{c(\sigma)}(x) w_\sigma(x) d\omega \le c \, 
    \frac{\|f\|_{\kappa,1}}{\sigma}.  
$$
\end{thm}

\begin{proof}
For each $x \in c(\sigma)$, there is a $\theta$ such that 
$$
 \int_{S^d} |f(y)| V_\kappa\left[ \chi_{c(x,\theta)}\right](y) 
    h_\kappa^2(y)d\omega \ge 
   \sigma \int_{S^d}  V_\kappa \left[\chi_{c(x,\theta)}\right](y)
   h_\kappa^2(y)d\omega.
$$
Assume that there is an $x$ for which $\theta \ge \pi/N$ for a fixed $N$,
say $N=20$. Then $\int_0^\theta (\sin \phi)^{2\lambda}d\phi \sim
  \int_0^\pi (\sin \phi)^{2\lambda}d\phi$.  
Since $V_\kappa \left[\chi_{c(x,\theta)}\right](x) \le V_\kappa 1 = 1$, it 
follows from \eqref{eq:2.8} that
\begin{align*}
 \int_{S^d} |f(y)| h_\kappa^2(y)d\omega
  & \ge  \int_{S^d} |f(y)| V_\kappa \left[\chi_{c(x,\theta)}\right](y)
    h_\kappa^2(y)d\omega \\
  & \ge \sigma \int_{S^d} V_\kappa \left[\chi_{c(x,\theta)}\right](y)
   h_\kappa^2(y)d\omega = \sigma \int_0^\theta (\sin \phi)^{2\lambda}d\phi\\
  & \ge c \sigma \int_0^\pi (\sin \phi)^{2\lambda}d\phi 
      \ge c \sigma \int_{S^d} d\omega \\
  & \ge c \sigma \int_{S^d}  \chi_{c(\sigma)}(y) d\omega,
\end{align*}
which is the stated inequality. Hence, we only need to consider the case that 
$0 \le \theta \le \pi/N$ for $N =20$. We use a covering lemma, which 
states that if $E$ is a subset of $S^d$ and $E$ is covered by a 
family of spherical caps $\{c_S(x,\theta)\}$, then a disjoint sequence 
$c_S(x_j,\theta_j)$, $j=1,2,\ldots$, can be chosen from the family such that 
$E \subset \cup_j c_S(x_j, 5 \theta_j)$. Such a lemma is proved in exactly the 
same way that the similar lemma with the solid ball in $\RR^d$ is proved 
(see \cite[p. 8]{St}). Let $c_S(x_j,\theta_j)$ be the sequence for the set 
$c(\sigma)$. Then
$$
  \sum_j \chi_{c_S(x_j,10 \theta_j)}(x) \ge c \chi_{c(\sigma)}(x), 
     \qquad x \in S^d,
$$
where we have enlarge the sets from $c_S(x_j,5\theta_j)$ to 
$c_S(x_j,10\theta_j)$, so that the inequality still holds and, furthermore,
the enlargement means that the solid set $c(x_j,10\theta_j) \supset
c(x_j,5\theta_j)$, from which it follows that the solid set 
$\cup_j c(x_j,10 \theta_j) \supset c^*(\sigma)$, where 
$$
c^*(\sigma):= \{r y': y' \in c(\sigma), 1-v(y') \le r\le 1\} 
$$ 
and $v(y') = \max \{r: r x' \in \cup_j c(x_j,10 \theta_j)\} >0$. The 
enlargement ensures that $v(y') > 0$ for all $y' \in c(\sigma)$. Thus,
$$
  \sum_j \chi_{c(x_j,10 \theta_j)}(x) \ge c \chi_{c^*(\sigma)}(x), 
     \qquad x \in B^{d+1}.
$$
By \eqref{eq:2.8}, it follows that for $\theta \le \pi /20$, 
\begin{align*}
 &\int_{S^d} V_\kappa \left[\chi_{c(x,\theta)}\right](y)h_\kappa^2(y)d\omega 
   = c \int_0^\theta (\sin\phi)^{2\lambda} d\phi \\
  & \qquad \ge c \int_0^{10 \theta} (\sin\phi)^{2\lambda} d\phi 
  =  c  \int_{S^d} V_\kappa \left[\chi_{c(x,10 \theta)}\right](y)
       h_\kappa^2(y)d\omega.
\end{align*}
Hence, by Fatou's lemma, it follows that 
\begin{align*}
\|f\|_{\kappa,1} & \ge \sum_j \int_{S^d} |f(y)| V_\kappa \left[
   \chi_{c(x_j,\theta_j)}\right](y) h_\kappa^2(y)d\omega \\ 
& \ge \sigma \sum_j \int_{S^d}V_\kappa \left[\chi_{c(x_j,\theta_j)}\right](y)
       h_\kappa^2(y)d\omega 
 \ge c \sigma  \int_{S^d} 
  \sum_j V_\kappa \left[\chi_{c(x_j,10 \theta_j)}\right](y)
       h_\kappa^2(y)d\omega.
\end{align*}
By the integral representation of $V_\kappa$ at \eqref{eq:2.5} and Fatou's 
lemma, $\sum_j V_\kappa \left[\chi_{c(x_j,10 \theta_j)}\right](y) \allowbreak
\ge V_\kappa \left[\sum_j \chi_{c(x_j,10 \theta_j)}\right](y)$. 
Consequently, using the inequality $\sum_j \chi_{c(x_j,10 \theta_j)}(x) 
\ge c \chi_{c^*(\sigma)}(x)$, we conclude that 
\begin{align*}
\|f\|_{\kappa,1} & \ge c\,\sigma \int_{S^d}  
  V_\kappa \left[\chi_{c^*(\sigma)}\right](y) h_\kappa^2(y)d\omega \\
& = c\,\sigma 
 \int_{B^{d+1}} \chi_{c^*(\sigma)}(x) (1-\|x\|^2)^{\gamma_\kappa -1}dx\\
& = c\,\sigma \int_0^1 r^{d-1} \int_{S^d} \chi_{c^*(\sigma)}(rx')d\omega(x')
    (1-r^2)^{\gamma_\kappa -1}dr \\ 
& = c\,\sigma \int_{S^d} \left[\int_0^1 r^{d-1} 
\chi_{c^*(\sigma)}(rx')(1-r^2)^{\gamma_\kappa -1}dr\right] d\omega(x'). 
\end{align*}
To complete the proof, we show that the inner integral is greater than
$w_\sigma(x') \chi_{c(\sigma)}(x')$ for a positive function $w_\sigma$ on 
$c(\sigma)$. By the definition of $c^*(\sigma)$, if $\chi_{c(\sigma)}(x')=1$,
then $\chi_{c^*(\sigma)}(rx') = 1$ for $1-v(x') \le r \le 1$, so that
\begin{align*}
 \int_0^1 r^{d-1} \chi_{c^*(\sigma)}(rx') (1-r^2)^{\gamma_\kappa -1}dr 
 & = \int_{1-v(x')}^1 r^{d-1} (1-r^2)^{\gamma_\kappa -1}dr \\
 & \ge c \int_{1- v(x')}^1 (1-r)^{\gamma_\kappa -1} dr 
   = c [v(x')]^{\gamma_\kappa}.
\end{align*}
Let $w_\sigma(x') = [v(x')]^{\gamma_\kappa}$. Then $w_\sigma(x') > 0$ for
all $x' \in c(\sigma)$. This completes the proof. 
\end{proof}

The inequality proved in the theorem is not the usual weak type $(1,1)$  
estimate because of the presence of $w_\sigma$ and the fact that $w_\sigma$
depends on $f$. Since $\|\CM_\kappa f\|_\infty \le \|f\|_\infty$, an ordinary
weak type inequality would imply, by the Marcinkiewicz interpolation theorem, 
that $\CM_\kappa f$ is of strong type $(p,p)$. 

\begin{cor} \label{cor:ae}
For $\theta \ge 0$ and $x \in S^d$, define
$$
  f_\theta(x) =  
  \frac{\int_{S^d} f(y) V_\kappa \left[\chi_{c(x,\theta)}\right](y)
    h_\kappa^2(y)d\omega} {\int_{S^d}
      V_\kappa \left[\chi_{c(x,\theta)}\right](y)
    h_\kappa^2(y)d\omega}.
$$
If $f \in L^1(h_\kappa^2)$, then $\lim_{\theta \to 0} f_\theta(x) = f(x)$ 
for almost every $x \in S^d$. 
\end{cor}

\begin{proof}
From the definition of $T_\theta^\kappa f$, an alternative definition of 
$f_\theta(x)$ is 
\begin{equation}\label{eq:2.10}
 f_\theta(x)  = \frac{\int_0^\theta T_\phi^\kappa f(x) 
   (\sin \phi)^{2\lambda} d\phi} {\int_0^\theta (\sin \phi)^{2\lambda} d\phi}.
\end{equation} 
Hence, it follows from 
$$
  f_\theta(x) - f(x) = \frac{\int_0^\theta (T_\phi^\kappa f(x) -f(x))
   (\sin \phi)^{2\lambda} d\phi}{\int_0^\theta (\sin \phi)^{2\lambda}d\phi} 
$$
and the inequality \eqref{eq:2.7} on $T_\theta^\kappa f$ that 
$$
  \|f_\theta - f\|_{\kappa,p} \le 
     \sup_{0 \le \phi \le \theta} \|T_\phi^\kappa f - f\|_{\kappa,p}
    :=\omega(f;\theta)_{\kappa,p}
$$
for $1 \le p \le \infty$. It is shown in \cite{X02b} that 
$\omega(f;\theta)_{\kappa,p} \to 0$ as $\theta \to 0$ for 
$f\in L^p(h_\kappa^2)$. In particular, it follows that $f_\theta \to f$ in 
$L^1(h_\kappa^2)$. The rest of the proof follows from the standard argument 
(see, for example, \cite[p. 8]{St}).  We need to show that for 
$f \in L^1(h_\kappa^2)$ and for almost every 
$x \in S^d$,
$$
  \Omega f(x) : = \left|\limsup_{\theta \to 0} f_\theta(x) - 
        \liminf_{\theta \to 0}f_\theta(x) \right| =0.
$$
If $g$ is continuous, $\Omega g(x) \equiv 0$. If $g \in L^1(h_\kappa^2)$, 
then 
$$
   \int_{S^d} \chi_{c(\sigma)}(y) w_\sigma(y) d\omega \le c 
    \frac{\|g\|_{\kappa,1}}{\sigma},
$$
where $c(\sigma) = \{x: 2 \CM_\kappa g(x) \ge \sigma\}$. Since $\Omega g(x)\le
2\CM_\kappa g(x)$ implies $\{x: \Omega g(x) \ge \sigma\} \subset c(\sigma)$,
it follows that 
$$
 \int_{S^d} \chi_{\{\Omega g(x) >\sigma\}}(y) w_\sigma(y) d\omega(y) 
 \le c  \frac{\|g\|_{\kappa,1}}{\sigma}.
$$
Every $f \in L^1(h_\kappa^2)$ can be written as a sum of $f = h + g$ with
$h \in C(S^d)$ and $\|g\|_{\kappa,1}$ arbitrarily small. Since $\Omega h(x)
\equiv 0$, the above inequality implies that $\chi_{\{\Omega f(x) >\sigma\}}
(y) w_\sigma(y)=0$ almost everywhere for all $\sigma >0$. Since $w_\sigma$ is 
positive on $c(\sigma)$, hence positive on $\{x:\Omega f(x) >\sigma\}
= \{x:\Omega g(x) >\sigma\} \subset c(\sigma)$, it 
follows that $\Omega f(x) = 0$ almost everywhere. 
\end{proof}

\subsection{Almost everywhere convergence of $h$-harmonic expansions} 

We will establish almost everywhere convergence for summation methods
of $h$-harmonic expansions. In \cite{X02b} a convolution, $\star_\kappa$, is 
defined for $f \in L^1(h_\kappa^2)$ and $g \in L^1(w_\lambda;[-1,1])$, 
\begin{equation} \label{eq:2.11} 
 (f\star_\kappa g)(x): = a_\kappa \int_{S^d} f(y) 
   V_\kappa[g(\langle x,\cdot\,\rangle)](y) h_\kappa^2(y) d\omega.
\end{equation} 
If $h_\kappa(x) \equiv 1$ (the usual surface measure), this is the spherical
convolution in \cite{CZ}. Because of \eqref{eq:2.4} and \eqref{eq:2.2}, a 
summation method for $h$-harmonic expansions can be written in the form of 
\begin{equation} \label{eq:2.12} 
Q_r f(x) = (f\star_\kappa q_r)(x), \qquad
   q_r(t) = \sum_{j=0}^\infty a_j(r)\frac{j+\lambda_\kappa}{\lambda_\kappa}
       C_j^{\lambda_\kappa}(t),
\end{equation} 
where $q_r$ is a function defined on $[-1,1]$, $r$ is a parameter, and $q_r$
satisfies  
$$ 
  a_\kappa \int_{S^d} V_\kappa q_r(\langle x, \cdot \rangle) 
  h_\kappa^2(y)d\omega = c_\lambda \int_0^\pi q_r(\cos\theta) 
   (\sin\theta)^{2\lambda} d\theta =1,  
$$
in which the first equal sign follows from \eqref{eq:2.8}. The summation method
$Q_r f(x)$ converges to $f(x)$ in $L^p(h_\kappa^2)$ for $1 \le p < \infty$  
and in $C(S^d)$ if $p = \infty$, as $r \to r_0$, where $r_0$ can be infinity. 

Our goal is to prove the almost everywhere convergence of $Q_r f(x)$. The 
following is a theorem that holds under mild assumptions on the kernel 
function.  

\begin{thm}\label{thm:a.e}
Assume that $|q_r(\cos\theta)| \le m_r(\theta)$ for some nonnegative 
differentiable function $m_r$ on $[0,\pi]$ and $m_r$ satisfies 
\begin{enumerate}
\item{}  $m_r'(\theta) \le 0$ for $0 \le \theta \le \pi$, 

\medskip

\item{} $\displaystyle{\sup_r c_\lambda \int_0^\pi m_r(\theta) 
 (\sin \theta)^{2\lambda}  d\theta < \infty}$, \quad $\lambda =\lambda_\kappa$.
\end{enumerate}
Then for $f \in L^1(h_\kappa^2)$ and all $x \in S^d$, 
$$
 \sup_r|Q_rf(x)| \le c \CM_\kappa f(x).
$$
In particular, $Q_r f(x)$ converges to $f(x)$ for almost every $x \in S^d$. 
\end{thm} 

\begin{proof}
Define
$$
   \Lambda(\theta, x) =   \int_0^\theta T_\phi^\kappa |f|(x) 
      (\sin \phi)^{2\lambda} d \phi. 
$$
It follows from the definition of $\CM_\kappa f(x)$ and \eqref{eq:2.10} that
$$
\Lambda(\theta,x) \le \CM_\kappa f(x) \int_0^\theta(\sin\phi)^{2\lambda}d\phi
$$
for all $x \in S^d$. By the definition of $T_\theta^\kappa f$ and the fact 
that $|T_\theta^\kappa f(x)|\le T_\theta^\kappa |f|(x)$, 
$$ 
 |Q_r f(x)| = c_\lambda \left | \int_0^\pi T_\phi^\kappa f(x) 
             g_r(\cos \phi) (\sin \phi)^{2\lambda} d \phi \right|
   \le c_\lambda \int_0^\pi T_\phi^\kappa |f|(x) 
             m_r(\phi) (\sin \phi)^{2\lambda} d \phi.
$$
Integrating by parts, we obtain
\begin{align*}
 |Q_r f (x)| & \le c_\lambda \left[\Lambda(\pi,x) m_r(\pi) -
          \int_0^\pi \Lambda(\theta,x) m_r'(\theta) d\theta \right]\\ 
      & \le c_\lambda \CM_\kappa f(x) \left[m_r(\pi)\int_0^\pi 
        (\sin \phi)^{2\lambda}d\phi
       -   \int_0^\pi m_r'(\theta) \int_0^\theta (\sin \phi)^{2\lambda} 
        d\phi d\theta \right], 
\end{align*}
since $m_r'(\cos \theta) \le 0$. Integrating by parts again, we conclude that 
$$
  |Q_r f(x)| \le \CM_\kappa f(x) c_\lambda \int_0^\pi m_r(\theta) 
     (\sin \theta)^{2\lambda} d\theta \le c  \CM_\kappa f(x),
$$
where $c$ is independent of $r$. This establishes the stated maximal 
inequality, from which the almost everywhere convergence of $Q_rf(x)$ 
follows as in the proof of Corollary \ref{cor:ae}.  
\end{proof}

\medskip\noindent
{\it Remark.} From the proof it follows that we can replace the assumption
on $m(\theta)$ by 
$$
\sup_r m_r(\pi) \le c \qquad\hbox{and}\qquad   \sup_r c_\lambda 
  \int_0^\pi |m_r'(\theta)| (\sin \theta)^{2\lambda}  d\theta \le c.
$$

We apply the above theorem to three summation methods for $h$-harmonic 
expansions. The first one is the analog of the Poisson integral defined by 
$$
 P_r(f; x) = (f \star_\kappa q_r)(x), \qquad
  q_r(\cos \theta) = \frac{1-r^2}{(1-2r \cos \theta + r^2)^{\lambda_\kappa+1}}.
$$ 
for $x \in S^d$ and $r < 1$. The kernel function of this integral is the 
Poisson kernel for the $h$-harmonics, which is equal to $\sum_{n=0}^\infty r^n 
P_n(h_\kappa^2;x,y)$. It is proved in \cite{X02b} that $P_r(f;x)$ converges
to $f(x)$ in $L^p(h_\kappa^2)$ as $r \to 1-$. 

\begin{prop} \label{prop:poisson}
For $f \in L^1(h_\kappa^2)$, 
$$
   \lim_{r \to 1-} P_r(f,x) = f(x), \qquad \hbox{a.e. on $S^d$.} 
$$
\end{prop}

\begin{proof}
Clearly $m_r(\theta)=q_r(\cos \theta) \ge 0$ and it is easy to see that 
$m_r'(\theta) \le 0$. Hence we can apply the theorem. In fact, in this case
we have $|P_r(f,x)| \le \CM_\kappa f(x)$.  
\end{proof}

Our second example is the de la Vall\`ee Poussin means defined by
\begin{equation} \label{eq:2.13}
M_n(f,x):= \sum_{k=0}^n \mu_{k,n}^\lambda Y_n(h_\kappa^2;f,x) = 
  (f\star_\kappa q_n)(x) 
\end{equation}
where 
$$
 \mu_{k,n}^\lambda = \frac{n!}{(n-k)!} 
          \frac{\Gamma(n+ 2 \lambda + 1)}{\Gamma(n+k+2 \lambda+1)}, \qquad
   \lambda = \lambda_\kappa.
$$ 
Using the compact formula for the reproducing kernel and the formula for 
the Gegenbauer polynomials (see, for example, \cite[p. 11]{Askey}), 
$$
 \sum_{k=0}^n \frac{n!}{(n-k)!} \frac{\Gamma(n+ \lambda+1/2)}
   {\Gamma(n+k+2 \lambda+1)} \frac{k+\lambda}{\lambda} C_k^\lambda(t) 
    = \frac{\Gamma(1/2)}{2^{2 \lambda} \Gamma(\lambda+1)} 
       \left(\frac{1+t}{2}\right)^n,
$$
the kernel function is given by 
$$
q_n(\cos \theta) = \frac{(2\lambda+1)_n}{(\lambda+1/2)_n} \left(
    \frac{1+\cos \theta}{2}\right)^n
   = \frac{(2\lambda+1)_n}{(\lambda+1/2)_n} 
     \left( \cos \frac{\theta}{2}\right)^{2n},
$$
where $(a)_n = a (a+1) \ldots (a +n-1)$. 
Clearly $m_n(\theta) = q_n(\cos \theta) \ge0$ and it is easy to see that 
$m_n'(\theta) \le 0$. Hence, we can apply Theorem \ref{thm:a.e} to conclude 
that $M_n f(x)$ converges almost everywhere.

\begin{prop} \label{Poussin}
For $f \in L^1(h_\kappa^2)$, the de la Vall\'ee Poussin means $M_n f$ 
satisfy
$$
   \lim_{n \to \infty} M_n(f,x) = f(x), \qquad \hbox{a.e. on $S^d$.} 
$$
\end{prop}

For the ordinary harmonic expansions, this was proved in \cite{BL}. These 
means were introduced by de la Vall\`ee Poussin for Fourier series, and they
have been extended to various other series, such as Gegenbauer series and 
Jacobi series. See \cite{Bav,BL,X02} for further references. 

Our third example is the Ces\`aro $(C,\delta)$ means. For $\delta >0$, the 
Ces\`aro $(C, \delta)$ means, $s_n^\delta$, of a sequence $\{c_n\}$ are 
defined by 
$$
s^\delta_n = (A_n^\delta)^{-1}\sum_{k=0}^n A_{n-k}^\delta c_k, 
  \qquad A_{n-k}^\delta = \binom{n-k+\delta}{n-k}. 
$$
We denote the $n$-th $(C,\delta)$ means of the $h$-harmonic expansion by
$S_n^\delta(h_\kappa^2;f)$. These means can be written as
$$
 S_n^\delta(h_\kappa^2;f) = (f\star_\kappa q_n^\delta)(x), \qquad
   q_n^\delta(t) = (A_n^\delta)^{-1}\sum_{k=0}^n A_{n-k}^\delta 
      \frac{(k+\lambda)}{\lambda} C_k^\lambda(t),
$$
where $\lambda = \lambda_\kappa$. The function $q_n^\delta(t)$ is the kernel
of the $(C,\delta)$ means of the Gegenbauer expansions at $x =1$.

In both of the above examples, the kernel functions, $q_r(t)$, are positive 
and we can apply Theorem \ref{thm:a.e} with $m_r (\theta)= q_r(\cos \theta)$. 
The Ces\`aro $(C,\delta)$ means are positive if $\delta \ge 2 \lambda +1$.
It is proved in \cite{X97b} that, for $h_\kappa$ associated with every 
reflection group, $S_n^\delta(h_\kappa^2;f)$ converges to $f$ in 
$L^p(h_\kappa^2)$ if $\delta > \gamma_\kappa + (d-1)/2 = \lambda_\kappa$. Here
we have

\begin{prop} 
For $f \in L^1(h_\kappa^2)$ and $\delta > \lambda_\kappa$, the Ces\`aro
$(C,\delta)$ means satisfy  
$$
   \lim_{n \to \infty} S_n^\delta(h_\kappa^2;f,x) = f(x), \qquad 
      \hbox{a.e. on $S^d$.} 
$$
\end{prop}

\begin{proof}
It is known that the kernel $q_n^\delta$ is bounded by
$$
 | q_n^\delta(\cos \theta)| \le c \begin{cases} 
        n^{\lambda -\delta}(n^{-2} + \theta^2)^{(-\lambda -\delta -1)/2}, 
           & \delta \le \lambda+1 \\
        n^{-1}(n^{-2}+\theta^2)^{-\lambda -1}, & \delta \ge \lambda+1
    \end{cases}
$$
(see, for example, \cite{BC}).  Let $m_n(\theta)$ be the function in the right 
hand side of the above estimate. It is easy to see that $m_n'(\theta)$ is 
non-positive for $0 \le \theta \le \pi$. Furthermore, it is not hard to see
that $\int_0^\pi m_n(\theta)(\sin \theta)^{2\lambda}d\theta$ is uniformly 
bounded. Consequently, we can apply Theorem \ref{thm:a.e} to finish the 
proof.
\end{proof}

Let us mention that the index $\lambda_\kappa$ may not be the critical index 
for the $(C,\delta)$ means. In fact, in the case of $h_\kappa$ in 
\eqref{eq:1.2} associated with $\ZZ_2^{d+1}$, it is proved in \cite{LX} 
that the $(C,\delta)$ means $S_n^\delta(h_\kappa^2;f)$ converges to $f$ in 
$L^p(h_\kappa^2)$ if and only if $\delta > (d-1)/2+ \sum \kappa_i  
- \min_i \kappa_i$; that is, the critical index in the case of $G=\ZZ_2^{d+1}$ 
is $\lambda_\kappa - \min_i \kappa_i$.  The above proposition shows that 
$S_n^\delta(h_\kappa^2;f)$ converges almost everywhere to $f$ if $\delta >
\lambda_\kappa$. Naturally, we expect that $\lambda_\kappa - \min_i \kappa_i$
is also the critical index for the almost everywhere convergence. However, 
our method does not seem to be enough to prove such a result. 

\section{Maximal function and almost everywhere convergence on $B^d$}
\setcounter{equation}{0}

Recall the weight function $W_{\kappa,\mu}^B(x)$ defined in \eqref{eq:1.4},
in which $h_\kappa$ is an reflection invariant weight function defined on 
$\RR^d$. Let $a_{\kappa,\mu}$ denote the normalization constant for 
$W_{\kappa,\mu}^B$. Denote by $L^p(W_{\kappa,\mu}^B)$, $1 \le p \le \infty$, 
the space of measurable functions defined on $B^d$ with the finite norm
$$
  \|f\|_{W_{\kappa,\mu}^B,p}:= \Big(a_{\kappa,\mu} \int_{B^d} |f(x)|^p 
     W_{\kappa,\mu}^B(x) dx \Big)^{1/p}, \qquad 1 \le p < \infty,
$$
and for $p = \infty$ we assume that $L^\infty$ is replaced by $C(B^d)$, the 
space of continuous function on $B^d$.

Let $\CV_n^d(W_{\kappa,\mu}^B)$ denote the space of orthogonal polynomials of 
degree $n$ with respect to $W_{\kappa,\mu}^B$ on $B^d$. Elements of 
$\CV_n^d(W_{\kappa,\mu}^B)$ are closely related to the $h$-harmonics associated
with the weight function 
$$
h_{\kappa,\mu}(y_1,\ldots,y_{d+1}) = h_\kappa(y_1,\ldots,y_d)|y_{d+1}|^{\mu}
$$
on $\RR^{d+1}$, where $h_\kappa$ is associated with the reflection group $G$. 
The function $h_{\kappa,\mu}$ is invariant under the group $G \times \ZZ_2$. 
Let $Y_n$ be an $h$-harmonic polynomial of degree $n$ associated to 
$h_{\kappa,\mu}$ and assume that $Y_n$ is even in its $(d+1)$-th variable; 
that is, $Y_n(x,x_{d+1})=Y_n(x,-x_{d+1})$. We can write
\begin{equation} \label{eq:3.1}
   Y_n(y) = r^n P_n(x), \qquad  y = r(x,x_{d+1}) \in \RR^{d+1}, \quad 
 r = \|y\|, \quad (x,x_{d+1}) \in S^d,
\end{equation}
in polar coordinates. Then $P_n$ is an element of $\CV_n^d(W_{\kappa,\mu}^B)$  
and this relation is an one-to-one correspondence \cite{X01a}. The intertwining
operator associated with $h_{\kappa,\mu}$, denoted by $V_{\kappa,\mu}$, 
is given in terms of the intertwining operator $V_\kappa$ associated to
$h_\kappa$ and the operator $V_\mu^{\ZZ_2}$ associated to 
$h_\mu(x) =|x_{d+1}|^\mu$, $x \in \RR^{d+1}$, which is given explicitly by
\eqref{eq:2.4} (setting $\kappa_{d+1}=\mu$ and $\kappa_i=0$ for $1\le i\le d$).
That is,
$$
  V_{\kappa,\mu} f(x,x_{d+1}) = c_\mu \int_{-1}^1 
         V_\kappa [f(\cdot,x_{d+1}t)](x) (1+t)(1-t^2)^{\mu-1}dt,
$$
where $x \in \RR^d$. Since polynomials in $\CV_n^d(W_{\kappa,\mu}^B)$ 
correspond to $h$-harmonics that are even in the last coordinates, we 
introduce a modified operator
\begin{align} \label{eq:3.2}
  V_{\kappa,\mu}^B f(x,x_{d+1}) & :=  \frac{1}{2}
   \left[V_{\kappa,\mu} f(x,x_{d+1})
     + V_{\kappa,\mu} f(x, -x_{d+1}) \right] \\
  & =  c_\mu \int_{-1}^1 
         V_\kappa [f(\cdot,x_{d+1}t)](x) (1-t^2)^{\mu-1}dt \notag
\end{align}
acting on functions defined on $\RR^{d+1}$. We use this operator to define a 
convolution structure,  $\star_{\kappa,\mu}^B$, for weight $W_{\kappa,\mu}^B$
on $B^d$.

\begin{defn} 
For $f\in L^1(W_{\kappa,\mu}^B)$ and $g \in L^1(w_{\lambda+\mu})$,
define
$$
  (f \star_{\kappa,\mu}^B g)(x) = a_{\kappa,\mu} \int_{B^d} f(x) 
     V_{\kappa,\mu}^Bg(\langle X, \cdot \rangle)(Y)W_{\kappa,\mu}^B(y)dy,
$$
where $\lambda = \lambda_\kappa$, $X = (x, \sqrt{1-\|x\|^2})$ and 
$Y = (y, \sqrt{1-\|y\|^2})$. 
\end{defn}

The properties of this convolution can be derived from the corresponding 
convolution $f \star_{\kappa,\mu} g$ defined in \eqref{eq:2.11} (with respect 
to $h_{\kappa,\mu}$ instead of $h_\kappa$). In fact, we have the following
proposition.

\begin{prop}  \label{prop:3.1}
For $f\in L^1(W_{\kappa,\mu}^B)$ and $g \in L^1(w_{\lambda+\mu})$, 
$$
(f \star_{\kappa,\mu}^B g)(x) = (F\star_{\kappa,\mu} g)(x, \sqrt{1-\|x\|^2}), 
\qquad F(x,x_{d+1}) = f(x).
$$
\end{prop}  

\begin{proof}
Using the elementary identity for $P$ defined on $S^d$, 
\begin{equation} \label{eq:3.3}
   \int_{S^d} P(y) d\omega = \int_{B^d} \left[ P(x,\sqrt{1-\|x\|^2}\,)+
     P(x,-\sqrt{1-\|x\|^2}\,) \right]\frac{dx}{\sqrt{1-\|x\|^2}},
\end{equation}
this is an easy consequence of the equation \eqref{eq:3.2}.
\end{proof}

In particular, the above relation and Proposition 2.2 in \cite{X02b} shows 
that $f\star_{\kappa,\mu}^B g$ satisfies Young's inequality:
 
\begin{prop}  \label{prop:3.2}
Let $p,q,r \ge 1$ and $p^{-1} = r^{-1}+q^{-1}-1$. For $f \in 
L^q(W_{\kappa,\mu}^B)$ and $g \in L^r(w_{\lambda+\mu};[-1,1])$, 
$$
\|f\star_{\kappa,\mu} g\|_{W_{\kappa,\mu}^B,p} \le \|f\|_{W_{\kappa,\mu}^B,q} 
   \|g\|_{w_{\lambda+\mu},r}, 
$$
where $\|\cdot\|_{w_{\lambda+\mu},r}$ denotes the $L^r(w_{\lambda+\mu},[-1,1])$
norm.  
\end{prop}

We now define the generalized translation operator on $B^d$, which is an 
analog of the spherical means. It is defined implicitly via the convolution. 

\begin{defn}\label{defn3}
For $f\in L^1(W_{\kappa,\mu}^B)$, the generalized translation operator 
$T_\theta(W_{\kappa,\mu}^B;f)$ is defined implicitly by 
\begin{align} \label{eq:3.4}
c_{\lambda+\mu} \int_0^\pi T_\theta(W_{\kappa,\mu}^B;f,x)g(\cos \theta) 
   (\sin \theta)^{2\lambda+2\mu}d\theta = (f\star_{\kappa,\mu}^B g)(x).
\end{align}
where $\lambda = \lambda_\kappa$, for every $g \in L^1(w_{\lambda+\mu})$. 
\end{defn}

For $F \in L^1(h_{\kappa,\mu}^2)$, denote by $T_\theta^{\kappa,\mu} F$ the 
weighted spherical means associated with the weight function $h_{\kappa,\mu}$
on $S^{d+1}$. 

\begin{prop} \label{prop:3.4}
For each $x \in B^d$ the operator $T_\theta(W_{\kappa,\mu}^B;f,x)$ is a 
uniquely determined $L^\infty(w_\lambda)$ function in $\theta$. Furthermore, 
define $F(x,x_{d+1}) = f(x)$; then
\begin{equation} \label{eq:3.5}
 T_\theta(W_{\kappa,\mu}^B;f,x) = T_\theta^{\kappa,\mu} F(X), \qquad
   X = (x,\sqrt{1-\|x\|^2}), \quad x \in B^d. 
\end{equation} 
\end{prop}

\begin{proof}
For $x \in B^d$ and $(x,x_{d+1}) \in S^d$, using \eqref{eq:3.2} and the fact
that $F(y,y_{d+1}) = f(y)$ and $h_{\kappa,\mu}^2(y,y_{d+1})$ are both even
in $y_{d+1}$, it follows from Proposition \ref{prop:3.1} and the definition of 
$T_\theta^{\kappa,\mu}$ in Definition \ref{defn1} that 
\begin{align*}
 (f \star_{\kappa,\mu}^B g)(x)  & 
= a_{\kappa,\mu} \int_{S^d} F(u) V_{\kappa,\mu}
      g(\langle X,\cdot \rangle)(u) h_{\kappa,\mu}^2(u) d\omega
\\
&  = c_{\lambda+\mu} \int_{0}^\pi T_{\theta}^{\kappa,\mu}
  F(x,x_{d+1}) g(\cos \theta) (\sin \theta)^{2\lambda+2\mu}d \theta. 
\end{align*}
Comparing with the equation \eqref{eq:3.4} gives \eqref{eq:3.5}. It follows
from \eqref{eq:3.5} and the property of $T_\theta^{\kappa,\mu}$ that 
$T_\theta(W_{\kappa,\mu}^B;f)$ is an $L^\infty(w_\lambda)$ function and 
it is uniquely defined.
\end{proof} 
 
The properties of $T_\theta(W_{\kappa,\mu}^B;f,x)$ can be derived from those
of $T_\theta^{\kappa,\mu}f$ proved in \cite{X02b}. In particular, it will
allows us to define a modulus of smoothness. Since it will not be used below,
we will not go into that direction here.

We now turn our attention to the main focus of this paper, the almost 
everywhere convergence. First we define the corresponding maximal function. 
   
\begin{defn}\label{defn4}
For $f\in L^1(W_{\kappa,\mu}^B)$, the maximal function $\CM_{\kappa,\mu}^B f$ 
is defined by 
\begin{equation*}
 \CM_{\kappa,\mu}^B f(x) = \sup_{0\le \theta \le \pi} 
    \frac{\int_0^\theta T_\phi(W_{\kappa,\mu}^B; |f|, x) 
    (\sin \phi)^{2\lambda_\kappa+2\mu} d\phi}
      {\int_0^\theta (\sin \phi)^{2\lambda_\kappa+2\mu} d\phi}.
\end{equation*} 
\end{defn}

Since $T_\phi(W_{\kappa,\mu}^B; f)$ is related to $T_\theta^{\kappa,\mu}F$, the
maximal function $\CM_{\kappa,\mu}^B f$ should be related to 
$\CM_{\kappa,\mu}f$ defined in Definition \ref{defn2} with respect to 
$h_{\kappa,\mu}$. 

\begin{prop}  \label{prop:3.7}
For $f \in L^1(W_{\kappa,\mu}^B)$ define $F(x,x_{d+1}) = f(x)$. Then  
\begin{equation} \label{eq:3.6}
  \CM_{\kappa,\mu}^B f (x) = \CM_{\kappa,\mu}F (X), \qquad 
 X =(x,\sqrt{1-\|x\|^2}).    
\end{equation}
Furthermore, define 
$$
e(x,\theta) = \{(y,y_{d+1}) \in B^{d+1}: 
 \langle x, y\rangle + \sqrt{1-\|x\|^2}\, y_{d+1} 
  \ge \cos \theta, \quad y_{d+1} \ge 0\}, 
$$
then an alternative formula for $\CM_{\kappa,\mu}^B f$ is 
\begin{equation*}
 \CM_{\kappa,\mu}^B f(x) = \sup_{0\le \theta \le \pi} 
   \frac{\int_{B^d}|f(y)| V_{\kappa,\mu}^B \left[\chi_{e(x,\theta)}\right](Y)
      W_{\kappa,\mu}^B(y) dy}  
        {\int_{B^d} V_{\kappa,\mu}^B \left[\chi_{e(x,\theta)}\right](Y) 
      W_{\kappa,\mu}^B(y) dy}.  
\end{equation*} 
\end{prop}

\begin{proof}
The first equation is a direct consequence of the Proposition \ref{prop:3.4}
and the definitions of the two maximal functions. To prove the second equation,
by \eqref{eq:3.6} and Proposition \ref{prop:2.3}, it suffices to show
$$
 \int_{S^d} |F(u)|V_{\kappa,\mu} \left[\chi_{c(X,\theta)}\right](u) 
  h_{\kappa,\mu}^2(u)d\omega(u) = 
 2 \int_{B^d} |f(y)|V_{\kappa,\mu}^B \left[\chi_{e(x,\theta)}\right](Y) 
  W_{\kappa,\mu}^B(y)dy, 
$$
where $c(X,\theta) = \{(y,y_{d+1}): \langle x, y\rangle + \sqrt{1-\|x\|^2}\, 
 y_{d+1} \ge \cos \theta\}$. By \eqref{eq:3.2} and \eqref{eq:3.3}, it suffices
to show 
$$
  V_{\kappa,\mu}^B \left[\chi_{c(X,\theta)}\right](y,\sqrt{1-\|y\|^2}) = 
    2 V_{\kappa,\mu}^B \left[\chi_{e(x,\theta)}\right](y,\sqrt{1-\|y\|^2}). 
$$ 
Let $e_{-}(x,\sigma) = \{(y,y_{d+1}): (y,-y_{d+1}) \in e(x,\sigma)\}$. By the 
definition of $V_{\kappa,\mu}^B$, it is easy to see that 
$$
V_{\kappa,\mu}^B\left[\chi_{e_{-}(x,\theta)}\right](y,y_{d+1}) = 
    V_{\kappa,\mu}^B \left[\chi_{e(x,\theta)}\right](y, - y_{d+1}).  
$$
Since $c(X,\theta) = e(x,\theta)\cup e_{-}(x,\theta)$ and $e(x,\theta)\cap 
e_{-}(x,\theta)$ has measure zero, the desired equation follows from the fact 
that $\CV_{\kappa,\mu}^B f(x,x_{d+1})$ is an even function in $x_{d+1}$. 
\end{proof}

It should be mentioned that one could define $M_{\kappa,\mu}^B f$ by 
\eqref{eq:3.6} without introducing $T_\theta(W_{\kappa,\mu}^B;f)$. We 
choose the longer way for the purpose of studying the situation on the
simplex in the next section, where $T_\theta(W_{\kappa,\mu}^B;f)$ is 
needed. 

We note that the upper hemisphere $S_{upper}^d = \{x: \|x\| =1, 
x_{d+1}\ge 0\}$ 
is isomorphic to the unit ball $B^d$, since points in $S_{upper}^d$ can be 
written as $(x,\sqrt{1-\|x\|^2})$ with $x \in B^d$. In particular, if $E$
is a set in $S_{upper}^d$, then it is isomorphic to the set
$c_B(\sigma) = \{y \in B^d:(y,y_{d+1} \in c(\sigma)\}$ in $B^d$. 
 
\begin{thm}  \label{thm:3.8}
Let $\sigma > 0$ and $f\in L^1(W_{\kappa,\mu}^B)$. Define 
$c_B(\sigma) = \{x \in B^d: \CM_{\kappa,\mu}^B f(x) \allowbreak
\ge \sigma\}$. Then there is a function $w_{\sigma}^B$ which is positive 
on $c_B(\sigma)$ such that
$$
\int_{B^d} \chi_{c_B(\sigma)}(x) w_\sigma(x) \frac{dx}{\sqrt{1-\|x\|^2}} 
 \le c  \frac{\|f\|_{W_{\kappa,\mu}^B,1}}{\sigma}. 
$$
\end{thm}

\begin{proof}
Let $F(x,x_{d+1}) = f(x)$ as before. By \eqref{eq:3.6}, 
$\CM_{\kappa,\mu}^B f(x)=\CM_{\kappa,\mu} F(X)$, where $X = 
(x,\sqrt{1-\|x\|^2}\,)$, and the proof of Proposition \ref{prop:3.7} also  
shows that $\CM_{\kappa,\mu}^B f(x)=\CM_{\kappa,\mu} F(X_-)$, where $X_- =  
(x,-\sqrt{1-\|x\|^2}\,)$. Let $c(\sigma) = \{(x,x_{d+1}):\CM_\kappa 
F(x,x_{d+1}) \ge \sigma\}$ for $(x,x_{d+1}) \in S^d$. It follows that 
$\chi_{c_B(\sigma)}(x)= \chi_{c(\sigma)}(X)$ and $\chi_{c_B(\sigma)}(x)=
\chi_{c(\sigma)}(-X)$. Let $w_\sigma$ be a function that is positive on 
$c(\sigma)$. We define $w_\sigma^B(x) = w_\sigma(X) +w_\sigma(X_-)$. 
Then $w_\sigma$ is positive on $c_B(\sigma)$. Hence, using \eqref{eq:3.3} 
and the inequality with respect to $h_{\kappa,\mu}^2$ in 
Theorem \ref{thm:maximal}, we get 
\begin{align*} 
\int_{B^d} \chi_{c_B(\sigma)}(x) w_\sigma^B(x) \frac{dx}{\sqrt{1-\|x\|^2}}  
& = \int_{S^d} \chi_{c_B(\sigma)}(x) w_\sigma(x,x_{d+1})d \omega(x,x_{d+1})\\ 
& = \int_{S^d} \chi_{c(\sigma)}(y) w_\sigma(y) d \omega(y)\\   
& \le c \, \sigma 
      \frac{\|f\|_{\kappa,\mu,1}}{\sigma} 
 = c \, \sigma  \frac{\|f\|_{W_{\kappa,\mu}^B,1}}{\sigma}, 
\end{align*} 
which completes the proof. 
\end{proof} 

Just as in the case of the sphere, the weak type inequality of the maximal
function $\CM_{\kappa,\mu}^Bf$ allows us to prove results on almost everywhere 
convergence for summation methods of orthogonal expansions with respect to
$W_{\kappa,\mu}^B$ on the unit ball.

For $f \in L^2(W_{\kappa,\mu}^B)$, its orthogonal expansion is given by
$$
 L^2(W_{\kappa,\mu}^B) = \sum_{n=0}^\infty\bigoplus \CV_n^d(W_{\kappa,\mu}^B)
 : \qquad f = \sum_{n=0}^\infty \proj_n^{\kappa,\mu} f, 
$$
where $\proj_n^{\kappa,\mu}:L^2(W_{\kappa,\mu}^B)\mapsto
\CV_n^d(W_{\kappa,\mu}^B)$ 
is the projection operator, which can be written as an integral
$$
  \proj_n^{\kappa,\mu} f(x) = a_{\kappa,\mu} \int_{B^d} f(y) 
      P_n(W_{\kappa,\mu}^B; x,y) W_{\kappa,\mu}^B(y)dy, 
$$
where $P_n(W_{\kappa,\mu}^B; x,y)$ is the reproducing kernel of 
$\CV_n^d(W_{\kappa,\mu}^B)$. It is known (\cite{X01a}) that
\begin{align*}
P_n(W_{\kappa,\mu}^B; x,y) = \frac{1}{2}
   \Big[
& Y_n\left(h_{\kappa,\mu}^2; (x,\sqrt{1-\|x\|^2}),(y,\sqrt{1-\|y\|^2})\right)\\
 & + Y_n\left(h_{\kappa,\mu}^2; (x,\sqrt{1-\|x\|^2}),(y,-\sqrt{1-\|y\|^2})
 \right) \Big],  \notag
\end{align*}
Because of \eqref{eq:2.2} and  \eqref{eq:3.2}, we have 
\begin{equation} \label{eq:3.7}
 P_n(W_{\kappa,\mu}^B;f,x) = \frac{n+\lambda_\kappa+\mu}{\lambda_\kappa+\mu}
  V_{\kappa,\mu}^B [C_n^{\lambda_\kappa+\mu} (\langle \cdot, Y \rangle )](X),
\end{equation}
where $X =(x,\sqrt{1-\|x\|^2})$ and $Y =(y,\sqrt{1-\|y\|^2})$. Hence, a  
summation method of the orthogonal expansions with respect to 
$W_{\kappa,\mu}^B$ can be written in the form of 
$$
 Q_r^B f(x) = (f\star_{\kappa,\mu}^B q_r)(x), \qquad    q_r(t) = 
 \sum_{j=0}^\infty a_j(r)\frac{j+\lambda_\kappa+\mu}{\lambda_\kappa+\mu} 
    C_j^{\lambda_\kappa +\mu}(t). 
$$
The almost everywhere convergence of $Q_r^B f(x)$ can be proved with the 
help of $\CM_{\kappa,\mu}^Bf(x)$, just as in the case of summability on 
the sphere. For example, there is an analog of Corollary \ref{cor:ae}. Indeed,
let 
$$
   f_\theta^B(x) = \frac{\int_0^\theta T_\phi(W_{\kappa,\mu}^B; f, x) 
    (\sin \phi)^{2\lambda_\kappa+2\mu} d\phi}
    {\int_0^\theta (\sin \phi)^{2\lambda_\kappa+2\mu} d\phi}, \qquad
   x \in B^d;
$$
then it is easy to see that $f_\theta^B(x) = F_\theta(X)$, where $F$ and $X$
are defined in Proposition \ref{prop:3.7}. Consequently, by Corollary 
\ref{cor:ae}, $f_\theta^B$ converges almost everywhere to $f$ on $B^d$ for
all $f \in L^1(W_{\kappa,\mu}^B)$. Moreover, we have the following analog 
of Theorem \ref{thm:a.e}.

\begin{thm}\label{thm:3.9}
Assume that $|q_r(\cos\theta)| \le m_r(\theta)$ for some nonnegative 
differentiable function $m_r$ on $[0,\pi]$ and $m_r$ satisfies 
\begin{enumerate}
\item{}  $m_r'(\theta) \le 0$ for $0 \le \theta \le \pi$, 

\medskip

\item{} $\displaystyle{\sup_r c_\lambda \int_0^\pi m_r(\theta) 
 (\sin \theta)^{2\lambda+2\mu}  d\theta < \infty}$, \quad $\lambda = 
 \lambda_\kappa$.
\end{enumerate}
Then for $f \in L^1(W_{\kappa,\mu}^B)$ and all $x \in B^d$, 
$$
 \sup_r|Q_r^Bf(x)| \le c \CM_{\kappa,\mu}^B f(x).
$$
In particular, $Q_r^B f(x)$ converges to $f(x)$ for almost every $x \in B^d$. 
\end{thm} 

\begin{proof}
Using \eqref{eq:3.2} and \eqref{eq:3.3}, it is easy to see that  
$Q_r^B f(x) = Q_r F(x)$, where $Q_r $ is as in \eqref{eq:2.12} with
$h_{\kappa,\mu}$ in place of $h_\kappa$. Hence, this is a corollary of 
Theorem \ref{thm:a.e}.
\end{proof}

We can apply this theorem to various summation methods. For example, we can 
consider $P_r^B(f;x):=\sum_{n=0}^\infty r^n P_n(W_{\kappa,\mu};f)$ as 
$r\mapsto 1-$, which can be written as a Poisson type integral,
$$
 P_r^B(f;x) = a_{\kappa,\mu} \int_{B^d} f(y) V_{\kappa,\mu}^B 
  \left [ \frac{1-r^2}{(1-r \langle\cdot, Y \rangle +r^2)^{\lambda+1}}\right]
  (X) W_{\kappa,\mu}^B(x)dx.
$$ 

\begin{prop} \label{prop:3.10}
For $f \in L^1(W_{\kappa,\mu}^B)$, the Poisson summation with respect to 
$W_{\kappa,\mu}^B$ satisfies 
$$
   \lim_{n \to \infty} P_r^B(f;x) = f(x), \qquad 
      \hbox{a.e. on $B^d$.} 
$$
\end{prop}

The de la Vall\`ee Poussin means and the Ces\`aro means are defined for 
orthogonal expansions with respect to $W_{\kappa,\mu}^B$ just like the case
of $h$-harmonics on the sphere. Our result shows that they converge almost
everywhere on $B^d$ ($\delta > \lambda_\kappa +\mu$ for the Ces\`aro 
$(C,\delta)$ means). For $d =1$ and $h_\kappa(x) =1$, the expansions with
respect to $W_{\kappa,\mu}^B$ become the expansions in the Gegenbauer 
polynomials. 

\section{Maximal function and almost everywhere convergence on $T^d$}
\setcounter{equation}{0}

Recall the weight function $W_{\kappa,\mu}^T(x)$ defined in \eqref{eq:1.5},
in which $h_\kappa$ is a reflection invariant weight function defined on 
$\RR^d$ and $h_\kappa$ is even in each of its variables. The latter condition
means that the group $G$, under which $h_\kappa$ is invariant, has the 
abilian group $\ZZ_2^d$ as a subgroup. 
The definition of $L^p(W_{\kappa,\mu}^T)$, $1 \le p \le \infty$, is similar
to the case of $W_{\kappa,\mu}^B$. Similarly, the notions such as the space 
of orthogonal polynomials $\CV_n^d(W_{\kappa,\mu}^T)$ and the reproducing
kernel $P_n(W_{\kappa,\mu}^T;x,y)$ are defined as in the case of 
$W_{\kappa,\mu}^B$. 
 
Elements of $\CV_n^d(W_{\kappa,\mu}^T)$ are closely related to the polynomials 
in $\CV_{2n}^d(W_{\kappa,\mu}^B)$. Let us denote by $\psi$ the mapping 
$$ 
\psi: (x_1,\ldots,x_d) \in T^d \mapsto (x_1^2,\ldots,x_d^2) \in B^d
$$ 
and define $(f\circ \psi)(x_1,\ldots,x_d) = f(x_1^2,\ldots,x_d^2)$. Then
\begin{equation}\label{eq:4.1}
 P_{2n}(x) = (R_n \circ \psi)(x), \qquad x \in B^d,
\end{equation} 
is a one-to-one mapping between $R_n \in \CV_n^d(W_{\kappa,\mu}^T)$ and  
$P_{2n} \in \CV_{2n}^d(W_{\kappa,\mu}^B; \ZZ_2^d)$, the subspace of polynomials
in $\CV_{2n}^d(W_{\kappa,\mu}^B)$ that are even in each of its variables 
(invariant under $\ZZ_2^d$). This relation implies, in particular, that 
the reproducing kernel $P_n(W_{\kappa,\mu}^T;x,y)$ of
$\CV_n^d(W_{\kappa,\mu}^T)$ satisfies (\cite{X01b})
\begin{equation}\label{eq:4.2}
  P_n(W_{\kappa,\mu}^T;x,y) = \frac{1}{2^d} \sum_{\varepsilon \in \ZZ_2^d}
   P_{2n}\left(W_{\kappa,\mu}^B; x^{1/2}, \varepsilon y^{1/2}\right),
\end{equation}
where $x^{1/2} = (\sqrt{x_1}, \ldots, \sqrt{x_d})$ and $\varepsilon u = 
(\varepsilon_1 u_1,\ldots \varepsilon_d u_d)$. Using the fact that
$$
\frac{2 n+\lambda}{\lambda} C_{2n}^\lambda(t) = 
   p_n^{(\lambda-1/2,-1/2)}(1)p_n^{(\lambda-1/2,-1/2)}(2t^2-1), 
$$
where $p_n^{(\alpha,\beta)}$ denotes the orthonormal Jacobi polynomial of 
degree $n$, it follows from \eqref{eq:2.2} and \eqref{eq:4.2} that 
\begin{align*}
 P_n(W_{\kappa,\mu}^T;x,y) =\,
   &  p_n^{(\lambda_\kappa+\mu-\frac12,-\frac12)}(1)\\
   &  \times \frac{1}{2^d}  \sum_{\varepsilon \in \ZZ_2^d} V_{\kappa,\mu}^B 
      [p_n^{(\lambda_\kappa+\mu-\frac12,-\frac12)} 
         (2 \langle \cdot, \varepsilon Y^{1/2} \rangle^2-1)](X^{1/2}),
\end{align*}
where $X^{1/2} =(\sqrt{x_1},\ldots, \sqrt{x_d},\sqrt{1-|x|})$ and 
$\varepsilon Y^{1/2} (\varepsilon_1 \sqrt{x_1},\ldots, \varepsilon_d
   \sqrt{x_d},\sqrt{1-|x|})$. 
This formula suggests the following definition of a useful operator, 
$V_{\kappa,\mu}^T$, acting on functions of $d+1$ variables, 
\begin{equation}\label{eq:4.3}
  V_{\kappa,\mu}^T F(x,x_{d+1}) = \frac{1}{2^d}
    \sum_{\varepsilon \in \ZZ_2^d}
           V_{\kappa,\mu}^B F (\varepsilon x, x_{d+1}).
\end{equation}     
We use this operator to define a convolution operator on $T^d$: 

\begin{defn}\label{defn:4.1}
For $f \in L^1(W_{\kappa,\mu}^T)$ and $g \in L^1(w_{\lambda+\mu};[-1,1])$, 
we define 
$$
(f \star_{\kappa,\mu}^T g)(x) = a_{\kappa,\mu} \int_{T^d} f(y) 
  V_{\kappa,\mu}^T 
   \left[g \left(2 \langle X^{1/2},\cdot \rangle^2-1\right)\right]
   (Y^{1/2}) W_{\kappa,\mu}^T(y) dy, 
$$
where $X^{1/2} = \left(\sqrt{x_1},\ldots,\sqrt{x_d},\sqrt{1-|x|}\right)$. 
\end{defn}

Recall that $|x| = x_1 + \ldots + x_d$. Evidently, $f \star_{\kappa,\mu}^T g$ 
is related to the convolution structure $f \star_{\kappa,\mu}^B g$ on $B^d$. 
In fact, we have the following:

\begin{prop}
For $f \in L^1(W_{\kappa,\mu}^T)$ and $g \in L^1(w_{\lambda+\mu};[-1,1])$, 
$$
  \left((f\star_{\kappa,\mu}^T g)\circ \psi \right)(x) = 
     \left((f\circ \psi)\star_{\kappa,\mu}^B g(2\{\cdot\}^2-1) \right)(x).
$$
\end{prop} 

\begin{proof}
Using the elementary integral 
\begin{equation}\label{eq:4.4}
  \int_{B^d} f(x_1^2,\ldots,x_d^2) dx = \int_{T^d} f(x_1,\ldots,x_d)
    \frac{dx}{\sqrt{x_1\cdots x_d}}
\end{equation}
it is easy to see that $f \in L^1(W_{\kappa,\mu}^T)$ implies 
$f \circ \psi \in  L^1(W_{\kappa,\mu}^B)$. Furthermore, the equation
\eqref{eq:4.4} implies that 
\begin{align*}
  & (f \star_{\kappa,\mu}^T g)(x_1^2,\ldots,x_d^2)  = 
     a_{\kappa,\mu} \int_{B^d} (f\circ \psi)(y)  
      V_{\kappa,\mu}^T  \left[g \left(2 \langle X,\cdot \rangle^2-1\right)
        \right](Y) W_{\kappa,\mu}^B(y) dy \\
& \hspace{0.8in} =  a_{\kappa,\mu} \int_{B^d} (f\circ \psi)(y)  
    \frac{1}{2^d} \sum_{\varepsilon \in \ZZ_2^d} 
    V_{\kappa,\mu}^B  \left[g \left(2 \langle X,\cdot \rangle^2-1\right)
        \right](\varepsilon Y) W_{\kappa,\mu}^B(y) dy,
\end{align*}
where $\varepsilon Y = (\varepsilon_1 y_1,\ldots \varepsilon_d y_d, 1-|y|)$. 
Since $W_{\kappa,\mu}^B(y)$ is even in each of its variables, changing 
variables $y_i \mapsto \varepsilon_i y_i$ shows that the summation can be 
removed from the above formula.  
\end{proof}

Next we define the generalized translation operator associated with 
$W_{\kappa,\mu}^T$ on the simplex. It is again defined implicitly.

\begin{defn}\label{defn5}
For $f\in L^1(W_{\kappa,\mu}^T)$, the generalized translation operator 
$T_\theta(W_{\kappa,\mu}^T;f)$ is defined implicitly by 
\begin{align} \label{eq:4.5}
 c_{\lambda+\mu} \int_0^\pi T_\theta(W_{\kappa,\mu}^T;f,x)g(\cos 2 \theta) 
   (\sin \theta)^{2\lambda+2\mu}d\theta = (f \star_{\kappa,\mu}^T g)(x),
\end{align}
where $\lambda = \lambda_\kappa$, for every $g \in L^1(w_{\lambda+\mu})$. 
\end{defn}

The following proposition shows that the generalized translation operator 
on $T^d$ is closely related to the one on $B^d$, and it also shows that the
operator $T_\theta(W_{\kappa,\mu}^T;f)$ is well-defined. 

\begin{prop} \label{prop:4.1}
For each $x \in T^d$, the operator $T_\theta(W_{\kappa,\mu}^T;f,x)$ is a 
uniquely determined $L^\infty$ function in $\theta$. Furthermore, 
\begin{equation} \label{eq:4.6}
 \left(T_\theta(W_{\kappa,\mu}^T;f)\circ \psi\right) (x) = 
    T_\theta(W_{\kappa,\mu}^B;f\circ \psi, x), \qquad  x \in T^d. 
\end{equation} 
\end{prop}

\begin{proof}
From the previous proposition, we get
\begin{align*}
 & \left((f\circ \psi)\star_{\kappa,\mu}^B g(2\{\cdot\}^2-1) \right)(x)  =
 (f\star_{\kappa,\mu}^T g)(x_1^2,\ldots,x_d^2) \\
 & \hspace{.5in} 
 = c_{\lambda+\mu} \int_0^\pi T_\theta(W_{\kappa,\mu}^T;f,x_1^2,\ldots,x_d^2)
   g(\cos 2 \theta) (\sin \theta)^{2\lambda+2\mu}d\theta.
\end{align*}
Since $\cos 2 \theta = 2 \cos^2\theta -1$, it follows that 
$$
\left( (f\circ \psi)\star_{\kappa,\mu}^B g\right)(x)  =
  c_{\lambda+\mu} \int_0^\pi T_\theta(W_{\kappa,\mu}^T;f,x_1^2,\ldots,x_d^2)
   g(\cos \theta) (\sin \theta)^{2\lambda+2\mu}d\theta.
$$
Comparing with the definition of $T_\theta(W_{\kappa,\mu}^T;f)$ in 
\eqref{eq:3.4}, the equation \eqref{eq:4.6} follows. This also shows that 
$T_\theta(W_{\kappa,\mu}^T;f)$ is uniquely defined.
\end{proof} 
 
We use the generalized translation operator to define the maximal function
associated to $W_{\kappa,\mu}^T$. 
   
\begin{defn}\label{defn6}
For $f\in L^1(W_{\kappa,\mu}^T)$, the maximal function $\CM_{\kappa,\mu}^T f$ 
is defined by 
\begin{equation*}
 \CM_{\kappa,\mu}^T f(x) = \sup_{0\le \theta \le \pi} 
    \frac{\int_0^\theta T_\phi(W_{\kappa,\mu}^T; |f|, x) 
    (\sin \phi)^{2\lambda_\kappa+2\mu} d\phi}
      {\int_0^\theta (\sin \phi)^{2\lambda_\kappa+2\mu} d\phi}.
\end{equation*} 
\end{defn}

This maximal function is closely related to the maximal function on $B^d$ 
defined in Definition \ref{defn4}. Recall that 
$Y^{1/2} = (\sqrt{y_1},\ldots,\sqrt{y_d},\sqrt{1-|y|})$ for $y \in T^d$. 

\begin{prop}  \label{prop:4.8}
For $f \in L^1(W_{\kappa,\mu}^T)$, 
\begin{equation} \label{eq:4.7}
(\CM_{\kappa,\mu}^T f)\circ \psi =   \CM_{\kappa,\mu}^B(f\circ \psi).
\end{equation} 
Furthermore, define a subset of $B^{d+1}$, 
$$
e^+(x,\theta) = \{(y,y_{d+1}): \langle x^{1/2}, y\rangle + \sqrt{1-|x|}\,  
  y_{d+1}  \ge \cos \theta, \,\, y_i \ge 0, \,\, 1\le i \le d+1\}, 
$$
then an alternative formula for $\CM_{\kappa,\mu}^T f$ is 
\begin{equation*}
 \CM_{\kappa,\mu}^T f(x) = \sup_{0\le \theta \le \pi} 
   \frac{\int_{T^d}|f(y)| V_{\kappa,\mu}^T 
   \left[\chi_{e^+(x,\theta)}\right](Y^{1/2})W_{\kappa,\mu}^T(y) dy}  
    {\int_{T^d} V_{\kappa,\mu}^T \left[ \chi_{e^+(x,\theta)}\right](Y^{1/2}) 
      W_{\kappa,\mu}^T(y) dy}  
\end{equation*} 
for $x \in T^d$.
\end{prop}

\begin{proof}
The equation \eqref{eq:4.7} is a simple consequence of \eqref{eq:4.6}.
Recall the definition of $e(x,\theta)$ in the Proposition \ref{prop:3.7}.
In order to prove the alternative formula for $\CM_{\kappa,\mu}^T f$, 
it is sufficient to establish that 
\begin{align*}
 & \int_{B^d} |f\circ \psi|(y) V_{\kappa,\mu}^B 
  \left[\chi_{e(x,\theta)}\right](Y)
    W_{\kappa,\mu}^B (y)dy \\
 & \hspace{.5in} = 
  2^d \int_{T^d} |f(y)| V_{\kappa,\mu}^T [\chi_{e^+(x^2,\theta)}](Y^{1/2})
    W_{\kappa,\mu}^T (y)dy,  
\end{align*}
where $Y = (y,\sqrt{1-\|y\|^2})$ and $x^2 = \psi(x) = (x_1^2, \ldots, x_d^2)$.
We first notice that it is easy to see
\begin{align*}
e(x,\theta)&= \bigcup_{\varepsilon \in \ZZ_2^d} e^+(\varepsilon x^2, \theta) \\
& = \bigcup_{\varepsilon \in \ZZ_2^d} \{(y,y_{d+1}): 
 \langle \varepsilon x, y\rangle + \sqrt{1-\|x\|^2}\, y_{d+1}\ge \cos \theta, 
   \,\, y_i \ge 0, \,\, 1\le i \le d+1\}, 
\end{align*}
and the intersection of two sets $e^+(\varepsilon x^2, \theta)$ and 
$e^+(\varepsilon' x^2, \theta)$, where $\varepsilon \ne \varepsilon'$, 
has measure zero. In particular, this implies that 
$$
\chi_{e(x,\theta)}(u) = \sum_{\epsilon \in \ZZ_2^d} 
    \chi_{e^+(\epsilon x^2,\theta)}(u) 
$$
for almost all $u$. Recall that $V_\kappa$ is associated with the group $G$
which has $\ZZ_2^d$ as a subgroup and that $V_\kappa$ satisfies 
$R(\sigma) V_{\kappa} =V_{\kappa}R(\sigma)$, where $R(\sigma)f(x):= 
f(\sigma x)$, for every $\sigma \in G$, which holds, in particular, for 
$\sigma = \varepsilon \in \ZZ_2^d$. Hence, using \eqref{eq:3.2}, we see that, 
\begin{equation} \label{eq:4.8}
  V_{\kappa,\mu}^B \left[\chi_{e^+(\varepsilon x^2,\theta)}\right](y, y_{d+1})
= V_{\kappa,\mu}^B \left[\chi_{e^+(x^2,\theta)}\right](\varepsilon y, y_{d+1}).
\end{equation} 
Consequently, we get 
\begin{align*}
&\int_{B^d} |(f\circ\psi)(y)| V_{\kappa,\mu}^B \left[\chi_{e(x,\theta)}
\right](Y)
W_{\kappa,\mu}^B (y) dy \\
& \hspace{.5in}
= \int_{B^d} |(f\circ\psi)(y)|  V_{\kappa,\mu}^B \left[
  \sum_{\varepsilon \in \ZZ_2^d}\chi_{e^+(\varepsilon x^2,\theta)}(Y) \right]
    W_{\kappa,\mu}^B (y) dy \\
& \hspace{.5in}
= \int_{B^d} |(f\circ\psi)(y)| \sum_{\varepsilon \in \ZZ_2^d}
   V_{\kappa,\mu}^B \left[\chi_{e^+(x^2,\theta)} \right](\varepsilon Y)
    W_{\kappa,\mu}^B (y) dy, 
\end{align*}
where $\varepsilon Y = (\varepsilon y, \sqrt{1-\|y\|^2})$. The sum in the 
last integral is even in each of its variables, which allows us to use 
\eqref{eq:4.4} to get that the above integral is equal to 
$$
 \int_{T^d} |(f\circ\psi)(y)| \sum_{\varepsilon \in \ZZ_2^d}
 V_{\kappa,\mu}^B \left[\chi_{e^+(x^2,\theta)}(\varepsilon Y^{1/2}) \right]
   W_{\kappa,\mu}^T (y) dy, 
$$
from which the desired result follows from the definition of 
$V_{\kappa,\mu}^T$ in \eqref{eq:4.3}. 
\end{proof}

\begin{thm}  \label{thm:4.9}
Let $\sigma > 0$ and $f\in L^1(W_{\kappa,\mu}^T)$. Define $c_T(\sigma) = 
\{x \in T^d: \CM_{\kappa,\mu}^T f(x) \allowbreak \ge \sigma\}$. Then there is
a function $w_{\sigma}^T$ which is positive on $c_T(\sigma)$ such that
$$
\int_{T^d} \chi_{e_T(\sigma)}(x) w_\sigma^T(x) 
   \frac{dx}{\sqrt{x_1\cdots x_d (1-|x|)}}  
   \le c  \frac{\|f\|_{W_{\kappa,\mu}^T,1}}{\sigma}. 
$$
\end{thm}

\begin{proof}
Let $c_B(\sigma) = \{x\in B^d: \CM_{\kappa,\mu}^B(f\circ \psi)(x)\ge\sigma\}$.
Using $(\CM_{\kappa,\mu}^T f)\circ \psi = \CM_{\kappa,\mu}^B(f\circ \psi)$, we
see that $\chi_{c^T(\sigma)}(x^2) = 
 \chi_{c^B(\sigma)}(x)$. Define $w_\sigma^B = w_\sigma^T \circ \psi$. 
Using \eqref{eq:4.4}, we get 
\begin{align*} 
 \int_{T^d} \chi_{c_T(\sigma)}(x) w_\sigma^T(x)
   \frac{dx}{\sqrt{x_1\cdots x_d (1-|x|)}} 
& = \int_{B^d} \chi_{c_T(\sigma)}(x^2) w_\sigma^T(x^2) 
      \frac{dx}{\sqrt{1-\|x\|^2}}   \\
& =\int_{B^d} \chi_{c_B(\sigma)}(x) w_\sigma^B(x)\frac{dx}{\sqrt{1-\|x\|^2}}\\
& \le c \frac{\|f\circ \psi\|_{W_{\kappa,\mu}^B,1}}{\sigma} =
   c \frac{\|f\|_{W_{\kappa,\mu}^T,1}}{\sigma}, 
\end{align*}
where we have used the inequality in Theorem \ref{thm:3.8}. 
\end{proof}


Again, the weak type inequality of the maximal function $\CM_{\kappa,\mu}^Tf$ 
can be used to prove almost everywhere convergence. In particular, an analog 
of Corollary \ref{cor:ae} holds for 
$$
  f_\theta^T(x) =  \frac{\int_{T^d} f(y) V_{\kappa,\mu}^T 
   \left[\chi_{e^+(x,\theta)}\right](Y^{1/2})W_{\kappa,\mu}^T(y) dy}  
    {\int_{T^d} V_{\kappa,\mu}^T \left[ \chi_{e^+(x,\theta)}\right](Y^{1/2}) 
      W_{\kappa,\mu}^T(y) dy};  
$$
that is, $f_\theta^T(x)$ converges almost everywhere to $f(x)$ on $T^d$ 
as $\theta \to 0$ for all $f \in L^1(W_{\kappa,\mu}^T)$. Indeed, the proof
of Proposition \ref{prop:4.8} shows that $f_\theta^T \circ \psi = 
(f\circ \psi)_\theta^B$, so that the result follows from the almost everywhere 
convergence of $f_\theta^B f$ on $B^d$. 

We denote also by $\proj_n^{\kappa,\mu}: L^2(W_{\kappa,\mu}^T)\mapsto 
\CV_n^d(W_{\kappa,\mu}^T)$ the orthogonal projection operator. For $f \in 
L^2(W_{\kappa,\mu}^T)$, it can be written as an integral
$$
  \proj_n^{\kappa,\mu} f(x) = a_{\kappa,\mu} \int_{T^d} f(y) 
      P_n(W_{\kappa,\mu}^T; x,y) W_{\kappa,\mu}^T(y)dy, 
$$
using the reproducing kernel of $\CV_n^d(W_{\kappa,\mu}^T)$ given in 
\eqref{eq:4.2}. Using the operator $V_{\kappa,\mu}^T$, the reproducing kernel
can be written as 
$$
 P_n(W_{\kappa,\mu}^T;x,y) = p_n^{(\lambda_\kappa+\mu-\frac12,-\frac12)}(1)
   V_{\kappa,\mu}^T  \left[p_n^{(\lambda_\kappa+\mu-\frac12,-\frac12)} 
         (2 \langle \cdot, Y^{1/2} \rangle^2-1)\right](X^{1/2}),
$$
where $X^{1/2} =(\sqrt{x_1},\ldots, \sqrt{x_d},\sqrt{1-|x|})$. The summation 
method of the orthogonal expansions with respect to $W_{\kappa,\mu}^T$ can be
written in the form of 
$$
 Q_r^T f(x) = (f\star_{\kappa,\mu}^T q_r)(x), \qquad x \in T^d,
$$
where $q_r(t)$ has an expansion in terms of the Jacobi polynomials
$$
g_r(t) = \sum_{j=0}^\infty a_j(r)
   p_j^{(\lambda_\kappa+\mu-\frac12,-\frac12)}(1)
    p_j^{(\lambda_\kappa+\mu-\frac12,-\frac12)}(t).  
$$
The almost everywhere convergence of $Q_r^T f(x)$ can be proved with the 
help of $\CM_{\kappa,\mu}^Tf(x)$, just as in the case of summability on 
the sphere. 

Let us point, however, that $Q_r^Tf(x)$ is not directly related to $Q_r^Bf(x)$,
since the kernel of $Q_r^T f$ is expanded into the Jacobi series, while 
the kernel of $Q_r^B f$ is expanded into the Gegenbauer series. The quadratic
transformation between $P_n^{(\lambda-1/2,-1/2)}(2t^2-1)$ and 
$C_{2n}^\lambda(t)$ does not preserve the summability. Consequently, the 
convergence of $Q_r^Tf(x)$ will not follow as a consequence of the 
convergence of $Q_r^B $. In fact, this is precisely the reason why we have 
taken the steps to define $T_\theta(W_{\kappa,\mu}^T;f)$. Since we have 
developed all necessary intermediate steps, an analog of Theorem \ref{thm:a.e} 
can be proved. 

\begin{thm}\label{thm:4.10}
Assume that $|q_r(\cos \theta)| \le m_r(\theta)$ for some nonnegative 
differentiable function $m_r$ on $[0,\pi]$ and $m_r$ satisfies 
\begin{enumerate}
\item{}  $m_r'(\theta) \le 0$ for $0 \le \theta \le \pi$, 

\medskip

\item{} $\displaystyle{\sup_r c_{\lambda+\mu} \int_0^\pi m_r(\theta) 
 (\sin \theta)^{2\lambda+2\mu}  d\theta < \infty}$.
\end{enumerate}
Then for $f \in L^1(W_{\kappa,\mu}^B)$ and all $x \in B^d$, 
$$
 \sup_r|Q_r^T f(x)| \le c \CM_{\kappa,\mu}^T f(x).
$$
In particular, $Q_r^T f(x)$ converges to $f(x) \in L^1(W_{\kappa,\mu}^T)$ for 
almost every $x \in B^d$. 
\end{thm} 

\begin{proof}
First we note that changing the variable $\theta \mapsto \pi - \theta$ in the 
\eqref{eq:4.5} shows that 
\begin{equation} \label{eq:4.9}
  T_{\pi - \theta}(W_{\kappa,\mu}^T; f,x) = T_\theta(W_{\kappa,\mu}^T; f,x)
\end{equation}  
since $T_\theta(W_{\kappa,\mu}^T; f,x)$ is uniquely determined. By the 
definition of $Q_r^T f$ and \eqref{eq:4.5},
\begin{align*}
 Q_r^T f(x) & = c_{\lambda+\mu} \int_0^\pi T_\phi(W_{\kappa,\mu}^T; f,x) 
             g_r(\cos 2\phi) (\sin \phi)^{2\lambda+2\mu} d \phi \\
& = c_{\lambda+\mu} \frac{1}{2}\int_0^{2\pi} T_{\phi/2}(W_{\kappa,\mu}^T; f,x) 
        g_r(\cos \phi) \left(\sin \frac{\phi}{2}\right)^{2\lambda+2\mu} d \phi.
\end{align*}
Split the integral over $[0,2\pi]$ into two integrals, one over $[0,\pi]$ and
the other over $[\pi,2\pi]$. Changing the variable $\phi \mapsto 2\pi - \phi$ 
in the second integral and using \eqref{eq:4.9} shows that the integral over
$[\pi,2\pi]$ is equal to the one over $[0,\pi]$. Consequently, we get
$$
 Q_r^T f(x) 
  = c_{\lambda +\mu}\int_0^{\pi} T_{\phi/2}(W_{\kappa,\mu}^T; f,x) 
        g_r(\cos \phi) \left(\sin \frac{\phi}{2}\right)^{2\lambda+2\mu} d \phi.
$$
Let us define 
$$
  \Lambda(\theta, x) = \int_0^{\theta} T_{\phi/2}(W_{\kappa,\mu}^T; |f|,x) 
    \left(\sin \frac{\phi}{2}\right)^{2\lambda+2\mu} d \phi.
$$
It follows from the definition of $M_{\kappa,\mu}^T f(x)$ (changing variable
$\phi \mapsto \phi/2$) that 
$$
  \Lambda(\theta,x) \le M_{\kappa,\mu}^T f(x)
     \int_0^\theta \left(\sin \frac{\phi}{2}\right)^{2\lambda+2\mu} d \phi.
$$
We can now follow the proof of Theorem \ref{thm:a.e} almost literally to 
finish the proof.
\end{proof}

Let us point out that for $d =1$ and $h_\kappa(x) = |x|^{\kappa}$ (the 
group $G = \ZZ_2$), the weight function $W_{\kappa,\mu}^T$ becomes the Jacobi 
weight function $t^{\kappa-1/2} (1-t)^{\mu-1/2}$ on $[0,1]$. Hence, our result
provides a way to prove results on almost everywhere convergence of the Jacobi 
series. 

We can apply Theorem \ref{thm:4.10} to various summation methods. For example, 
we can consider the Poisson sum (integral) 
$$
 P_r^T(f;x) :=\sum_{n=0}^\infty r^n P_n(W_{\kappa,\mu}^T;f) 
    =\left( f \star_{\kappa,\mu}^T P(r;\{\cdot\})\right)(x), \qquad x \in T^d,
$$ 
as $r\mapsto 1-$, where the Poisson kernel $P(r;t)$ of the Jacobi series
$$
  P(r;t) = \sum_{n=0}^\infty r^n p_n^{(\lambda_\kappa+\mu-\frac12,-\frac12)}(1)
    p_n^{(\lambda_\kappa+\mu-\frac12,-\frac12)}(t)  
$$
has an explicit representation in terms of the hypergeometric function
(\cite[p. 102, Ex. 19]{Bailey} and use \cite[Vol. 1, p. 64, 2.1.4(23)]{Er})
$$
  P(r;t) = \frac{(1-r)(1+r)^{\lambda+\mu}}{(1-2r t +r^2)^{\lambda+\mu+1}}
     {}_2F_1\left( \begin{matrix} \frac{\lambda+\mu}{2},
     \frac{\lambda+\mu-1}{2}\\
              \frac{1}{2}  \end{matrix}; \frac{2r (1+t)}{(1+r)^2} \right).
$$
This kernel is known to be positive. 

\begin{prop}
For $f \in L^1(W_{\kappa,\mu}^T)$, 
$$
   \lim_{r \to 1-} P_r^T(f;x) = f(x)  \qquad \hbox{a.e. on $T^d$}.
$$
\end{prop}

\begin{proof}
We apply the theorem with 
$$ 
  P(r;\cos\theta) \le c \frac{(1-r^2)}{(1-2r \cos \theta +r^2)^{\lambda+\mu+1}}
   := m_r(\theta),
$$
where the inequality follows from the explicit formula of $P(r;t)$. It is easy 
to verify that $m_r$ satisfies the conditions required in the theorem. 
\end{proof}

The Ces\`aro means for the orthogonal expansion on $T^d$ are defined just
as the means on $B^d$. We have the similar result that states: for $f \in 
L^1(W_{\kappa,\mu}^T)$ and $\delta > \lambda_\kappa + \mu$, the Ces\`aro 
$(C,\delta)$ means $S_n^\delta(W_{\kappa,\mu}^T;f,x)$ converge almost 
everywhere to $f(x)$ on $T^d$. 

The de la Vall\`ee Poussin mean, $M_n^T(f,x)$, on the simplex is defined by 
$$
M_n^T(f,x):= \sum_{k=0}^n \mu_{k,n}^{\lambda_\kappa+\mu} 
     P_n(W_{\kappa,\mu}^T;f,x),
$$  
where $\mu_{k,n}^{\lambda+\mu}$ is as given in \eqref{eq:2.13}. Just as in the
case of Proposition \ref{Poussin}, our result shows the following:

\begin{prop} \label{prop:4.11}
For $f \in L^1(W_{\kappa,\mu}^T)$, the de la Vall\`ee Poussin means satisfy
$$
   \lim_{n \to \infty} M_n(W_{\kappa,\mu}^T;f,x) = f(x), \qquad 
      \hbox{a.e. on $T^d$.} 
$$
\end{prop}

We state the result for these means since they have an interesting connection 
in the case of the classical weight function \eqref{eq:1.7} on the simplex. In 
terms of the Bernstein basis 
$$
  B_\alpha(x) := \binom{n}{\alpha} x^\alpha (1-|x|)^{n-|\alpha|}, \qquad
    \alpha \in \NN_0^d, \;\; |\alpha| = n, 
$$
on the simplex, the de la Vall\`ee Poussin means with respect to the weight
function \eqref{eq:1.7} can be expressed as 
$$
  M_n^T(f;x) = \sum_{|\alpha|\le n} 
  \frac{\int_{T^d} f(y) B_\alpha(y)W_{\kappa,\mu}^T dy}{
     \int_{T^d} B_\alpha(y)W_{\kappa,\mu}^T dy} B_\alpha(x)
$$
and they are known as the Bernstein-Durrmeyer operators. These operators 
were studied by several authors (see, for example, \cite{BSX,De,Di}); in the 
case of $d=1$, the Jacobi series, they were studied in \cite{Bav,BX}. As far 
as we know, Proposition \ref{prop:4.11} appears to be new even for $d =1$.

\bigskip
{\it Acknowledgement.} The author thanks a referee for his careful
review that pinpointed a sticky point in the manuscript; understanding the
point led to the weak type inequality in its present form.


\begin{thebibliography}{99} 

\bibitem{Askey}
        R. Askey,
        \textit{Orthogonal polynomials and special functions}, 
	Regional Conference Series in Applied Mathematics \textbf{21}, 
        SIAM, Philadelphia, 1975.

\bibitem{AW}
        R. Askey and S. Wainger,
        A convolution structure for Jacobi series, 
        \textit{Amer. J. Math.} \textbf{91} (1969), p. 463-485.

\bibitem{Bailey}
        W. N. Bailey, 
        \textit{Generalized hypergeometric series},
        Cambridge Univ. Press, Cambridge, UK, 1935.

\bibitem{Bav} 
        H. Bavinck
        \textit{Jacobi series and approximation},
        Mathematical Centre Tracts, No. 39. 
        Mathematisch Centrum, Amsterdam, 1972.

\bibitem{BBP} 
        H. Berens, P. L. Butzer and S. Pawelke, 
        Limitierungsverfahren von Reihen mehrdimensionaler Kugelfunktionen
        und deren Saturationsverhalten, 
        \textit{Publ. Res. Inst. Math. Sci. Ser. A.} \textbf{4} (1968), 
        201-268.

\bibitem{BL}
        H. Berens and Luoqing Li,         
        On the de la Vall\'ee Poussin means on the sphere,
        \textit{Results in Math.}, \textbf{24} (1993), 12-26.         

\bibitem{BSX}
        H. Berens, H. J. Schmid, and Yuan Xu,
        Bernstein-Durrmeyer polynomials on a simplex,
        \textit{J. Approx. Theory}, \textbf{68} (1992), 247-261.

\bibitem{BX}
        H. Berens and Yuan Xu,
        On Bernstein-Durrmeyer polynomials with Jacobi weights,
        in \textit{Approximation Theory and Functional Analysis}, 
        25-46, Aacdemic Press, New York, 1991.

\bibitem{BC}
        A. Bonami and J-L. Clerc,
        Sommes de Ces\`aro et multiplicateurs des d\'eveloppe-\allowbreak ments
        en harmoniques sph\'eriques,
        \textit{Trans. Amer. Math. Soc.} \textbf{183} (1973), 223--263.

\bibitem{CZ} 
        A. P. Calderon and A. Zygmund,
        On a problem of Mihlin,
        \textit{Trans. Amer. Math. Soc.}, \textbf{78} (1955), 209-224. 

\bibitem{De}
        M. M. Derriennic,
        On multivariate approximation by Bernstein-type polynomials, 
        \textit{J. Approx. Theory}, \textbf{45} (1985), 155-166.

\bibitem{Di}
        Z. Ditzian, 
        Multidimensional Jacobi-type Bernstein-Durrmeyer operators,
        \textit{Acta Sci. Math. (Szeged)}, \textbf{60} (1995), 225-243.       

\bibitem{D1} 
        C. F. Dunkl,  
	Differential-difference operators associated to reflection groups,
        \textit{Trans. Amer. Math. Soc.} \textbf{311} (1989), 167--183.

\bibitem{D2} 
	C. Dunkl,
	Integral kernels with reflection group invariance,
 	\textit{Canad. J. Math.} \textbf{43} (1991), 1213-1227.

\bibitem{DX}
        C. F. Dunkl and Yuan Xu,
        \textit{Orthogonal polynomials of several variables},
        Cambridge Univ. Press, 2001. 
 
\bibitem{Er}
	A. Erd\'elyi, W. Magnus, F. Oberhettinger, and F. G. Tricomi,
	\textit{Higher transcendental functions}, 
	McGraw-Hill, New York, 1953.

\bibitem{LW}
        Luoqing Li and Kunyang Wang,
        \textit{Harmonic analysis and approximation on the unit sphere}
        Science Press, Beijing, 2000.
 
\bibitem{LX}
        Zh.-K, Li and Yuan Xu,
        Summability of orthogonal expansions of several variables, 
        \- {\it J. Approx. Theory}, \textbf{122} (2003), 267-333.

\bibitem{LN}
        P. I. Lizorkin and S. M. Nikolskii,
        Approximation theory on the sphere,
        \textit{Proc. Steklov Inst. Math.}, \textbf{172} (1987),       
        295-302.         
 
\bibitem{P}
        S. Pawelke, 
        \"Uber Approximationsordnung bei Kugelfunktionen und algebraischen
        Polynomen, 
        \textit{T\^ohoku Math. J.}, \textbf{24} (1972), 473-486.

\bibitem{Ros} 
	M. R\"osler,
	Positivity of Dunkl's intertwining operator,
	\textit{Duke Math. J.\/},  \textbf{98} (1999), 445--463.

\bibitem{Rus}
        Kh. Rustamov, 
        On approximation of functions on the sphere,  
        \textit{Russian Acad. Sci. Izv. Math.}, \textbf{43} (1994), 311-329.

\bibitem{St} 
 	E. M. Stein,
 	\textit{Singular integrals and Differentiability properties of
        functions}, 
        Princeton Univ. Press, Princeton, NJ, 1971.

\bibitem{Szego}
	G. Szeg\"{o},
	\textit{Orthogonal Polynomials},  
	Amer. Math. Soc. Colloq. Publ. Vol.23, Providence, 4th edition,
        1975.


\bibitem{V}
	N. J. Vilenkin,
	\textit{Special functions and the theory of group representations},
        American Mathematical Society Translation of Mathematics Monographs 
        \textbf{22}, 	
        American Mathematical Society, Providence, RI, 1968.

\bibitem{X97b}
	Yuan Xu,
	Integration of the intertwining operator for $h$-harmonic polynomials
	associated to reflection groups, 
	\textit{Proc. Amer. Math. Soc.} \textbf{125} (1997), 2963--2973.

\bibitem{X01a}
	Yuan Xu,            
	Orthogonal polynomials and summability in Fourier orthogonal 
        series on spheres and on balls,  
        \textit{Math. Proc. Cambridge Phil. Soc.}, \textbf{31} (2001), 
        139-155.    

\bibitem{X01b} 
	Yuan Xu,            
	Generalized classical orthogonal polynomials on the ball and on 
	the simplex,  
        \textit{Constr. Approx.}, \textbf{17} (2001), 383-412.

\bibitem{X02} 
	Yuan Xu,            
        Approximation by means of $h$-harmonic polynomials on the unit sphere,
        \textit{Adv. in Comp. Math}, to appear.
       \texttt{http://www.math.uoregon.edu/\~{}yuan} (ApproxSph.ps.gz). 

\bibitem{X02b} 
	Yuan Xu,            
        Weighted approximation of functions on the unit sphere,
        \textit{Const. Approx.}, to appear.
       \texttt{http://www.math.uoregon.edu/\~{}yuan} (BestApp.ps.gz). 

\end{thebibliography}
\enddocument